\documentclass[11pt]{amsart}
\title[Indecomposability properties]{Indecomposability of free group
factors over nonprime subfactors and abelian subalgebras}
\author[M. B. \c Stefan]{Marius B. \c Stefan}
\address{UCLA Mathematics Department, Los Angeles, CA 90095-1555}
\email{stefan@math.ucla.edu}
\newtheorem{thm}{Theorem}[section]
\newtheorem{prop}[thm]{Proposition}
\newtheorem{lem}[thm]{Lemma}
\newtheorem{cor}[thm]{Corollary}
\theoremstyle{definition}
\newtheorem{defn}[thm]{Definition}

\theoremstyle{remark}
\newtheorem{rem}[thm]{Remark}
\newcommand{\mb}[1]{\mathbb{#1}}
\newcommand{\mc}[1]{\mathcal{#1}}
\DeclareMathOperator{\s}{sp} \DeclareMathOperator{\dist}{dist}
\date{}
\begin{document}
\begin{abstract}
We use the free entropy defined by D. Voiculescu to prove that the free
group factors can not be decomposed as closed linear spans of noncommutative
monomials in elements of nonprime subfactors or abelian $*$-subalgebras,
if the degrees of monomials have an upper bound depending on the number of
generators. The resulting estimates for the hyperfinite and abelian dimensions
of free group factors settle in the affirmative a conjecture of L. Ge and
S. Popa (for infinitely many generators).
\end{abstract}
\maketitle
\section{Introduction}
L. Ge and S. Popa defined (\cite{5}) for a given type
$\mbox{\!I\!I}_1$-factor $\mc{M}$ the following two quantities:
$\ell_h(\mc{M})=\mbox{min}\{f\in
\mb{N}\,\,|\,\,\exists\,\,\mbox{hyperfinite}\,\, \mc{R}_1,\ldots
,$ $\mc{R}_f\subset \mc{M}\,\,\mbox{such that
}\,\,\overline{\mbox{sp}}^w \mc{R}_1\mc{R}_2\ldots
\mc{R}_f=\mc{M}\}$, $\ell_a(\mc{M})=\mbox{min}\{f\in
{\mb{N}}\,\,|\,\,\exists\,\,$ $ \mbox{abelian}\,\,\mc{A}_1,$ $
\ldots ,\mc{A}_f \subset \mc{M}\,\,\mbox{such that}\,\,
\overline{\mbox{sp}}^w \mc{A}_1\mc{A}_2\ldots \mc{A}_f=\mc{M}\}$
(the $\mbox{min}$ considered is $\infty$ if $\mc{M}$ can not be
generated as stated) and conjectured that
$\ell_h(\mc{L}(\mb{F}_n))=\ell_a(\mc{L}(\mb{ F}_n))=\infty$ for
$n\geq 2$, where $\mc{L}(\mb{F}_n)$ is the type
$\mbox{\!I\!I}_1$-factor associated to the free group with $n$
generators.

We use the concept of free entropy introduced by D. Voiculescu in
his breakthrough paper \cite{14} to prove that the conjecture
mentioned above is true at least partially (for $n=\infty$) that
is, $\ell_h(\mc{}\mc{L}(\mb{F}_n)),\ell_a(\mc{L}(\mb{F}_n))\geq
[\frac{n-2}{2}]+1$ for all $4\leq n\leq\infty$. Actually, our
result is more general and it states that the free group factor
with $n$ generators can not be asymptotically generated
(Definitions \ref{hde} and \ref{ade}) as $$\lim_{\omega\rightarrow
0}\,^{||\cdot ||_2}\sum_{\overset{1\leq j_1,\ldots ,j_{t+1}\leq
f}{1\leq t\leq d}} \mc{N}_{j_1}^\omega \mc{Z}^\omega
\mc{N}_{j_2}^\omega \mc{Z}^\omega\ldots \mc{N}_{j_t}^\omega
\mc{Z}^\omega \mc{N}_{j_{t+1}}^\omega$$ or
$$\lim_{\omega\rightarrow 0}\,^{||\cdot ||_2}\sum_{\overset{1\leq
j_1,\ldots ,j_{t+1}\leq f}{1\leq t\leq d}} \mc{A}_{j_1}^\omega
\mc{Z}^\omega \mc{A}_{j_2}^\omega \mc{Z}^\omega\ldots
\mc{A}_{j_t}^\omega \mc{Z}^\omega \mc{A}_{j_{t+1}}^\omega$$ if
$\{\mc{N}_1^\omega,\ldots ,\mc{N}_f^\omega\}_\omega$ are nonprime
subfactors, $\{\mc{A}_1^\omega,\ldots ,\mc{A}_f^\omega\}_\omega$
are abelian $*$-subalgebras, $\{\mc{Z}^\omega\subset
\mc{L}(\mb{F}_n)\}_\omega$ are subsets containing $p$ self-adjoint
elements, and $f,d\geq 1$ are integers such that $n\geq p+2f+1$.
Note that $\mc{L}(\mb{F}_n)$ admits decompositions of this sort if
we allow $d=\infty$, for  example if
$\mc{Z}^\omega=\mc{Z}=\{1\},f=n$, $\mc{N}_1^\omega
=\mc{N}_1,\ldots ,\mc{N}_n^\omega =\mc{N}_n$ are $n$ distinct
copies of the hyperfinite type $\mbox{\!I\!I}_1$-factor $\mc{R}$
and $\mc{A}_1^\omega =\mc{A}_1,\ldots ,\mc{A}_n^\omega =\mc{A}_n$
are $n$ distinct copies of $L^{\infty}([0,1])$ (since
$\mc{L}(\mb{F}_n)$ is both the free product of $n$ copies of
$\mc{R}$ and the free product of $n$ copies of
$L^{\infty}([0,1])$, see \cite{16}). Note also that the
indecomposability of $\mc{L}(\mb{ F}_n)$ as
$\overline{\mbox{sp}}^w\mc{N}\mc{Z}\mc{N}$ implies the primeness
of its subfactors (\cite{11}). Indeed, according to V. Jones
(\cite{7}), if $\mc{N}$ is a subfactor of finite index in $\mc{M}$
then $\mc{M}$ decomposes as $\mc{N}e\mc{N}$ where $e$ is the Jones
projection. In particular, the indecomposability properties of
$\mc{L}(\mb{F}_n)$ over nonprime subfactors and abelian
subalgebras are preserved to its subfactors of finite index.
Recall that the Haagerup approximation property (\cite{54}) is
another property preserved to the free group subfactors. A first
example of a prime $\mbox{\!I\!I}_1$-factor (with a nonseparable
predual, though) was given by S. Popa (\cite{53}) and then L. Ge
proved (with a free entropy estimate) that the free group factor
$\mc{L}(\mb{F}_n)$ is prime $\forall 2\leq n<\infty$ (\cite{4}),
thus answering a question from \cite{19}.

Our results are based on estimates of free entropy that is,
estimates of volumes of various sets of matrix approximants
(matricial microstates). The paper has four parts. After
introduction, we prove the first estimate of free entropy and
reobtain then a result of D. Voiculescu (\cite{14}): if a free
family of $m$ self-adjoint noncommutative random variables can be
generated by noncommutative power series by another family of $n$
self-adjoint noncommutative random variables, then $n\geq m$
(Theorem \ref{prop2}). However, we show that the assumption of
freeness from \cite{14} is not essential and it can be dropped. As
a consequence, the number of self-adjoint generators with finite
entropy, which generate a $*$-algebra $\mc{A}$ {\it
algebraically}, is constant. In the third part we prove the
indecomposability of $\mc{L}(\mb{F}_n)$ (and of its subfactors of
finite index) over nonprime subfactors (Theorem \ref{prop4}) and
in the last section, the indecomposability over abelian
subalgebras (Theorem \ref{prop6}).

We give next a short account on Voiculescu's free probability
theory (\cite{13}, \cite{16}) and on his original concept of free
entropy (\cite{14}, \cite{15}). A type $\mbox{\!I\!I}_1$-factor
$\mc{M}$ endowed with its unique normalized, faithful, normal
trace $\tau$ is sometimes called a $W^*$-probability space. The
trace $\tau$ determines the $2$-norm on $\mc{M}$, $||x||_2=\tau
(x^*x)^\frac{1}{2}$ and the completion of $\mc{M}$ w.r.t. $||\cdot
||_2$ is denoted $L^2(\mc{M},\tau )$. An element $x\in \mc{M}$ is
a semicircular element if it is self-adjoint and if its
distribution is given by the semicircle law: $$\tau
(x^k)=\frac{2}{\pi}\int_{-1}^{1}t^k\sqrt{1-t^2}dt\,\,\forall k\in
\mb{N}.$$ A family $(\mc{A}_i)_{i\in I}$ of unital $*$-subalgebras
of $\mc{M}$ is a free family provided that $\tau (x_1x_2\ldots
x_n)=0$ whenever $\tau (x_k)=0$, $x_k\in \mc{A}_{i_k}\,\,\forall
1\leq k\leq n,\,\,i_1,\ldots ,i_n\in I$ and $i_1\not=
i_2\not=\ldots\not= i_n,\,\,n\in \mb{N}$. A set $\{x_i\}_{i\in
I}\subset \mc{M}$ is free if the family $(*$-alg$\{1,x_i\})_{i\in
I}$ is free. A free set $\{x_i\}_{i\in I}\subset \mc{M}$
consisting of semicircular elements is called a semicircular
system. If $\mb{F}_n$ is the free group with $n$ generators
($2\leq n\leq \infty$) then $\mc{L}(\mb{ F}_n)$ denotes
(\cite{32}) the von Neumann algebra generated by the left regular
representation $\lambda
:\mb{F}_n\rightarrow\mc{B}(l^2(\mb{F}_n))$. $\mc{L}(\mb{ F}_n)$
is a factor of type $\mbox{\!I\!I}_1$ - the free group factor on
$n$ generators. It has a canonical trace $\tau (\cdot )=(\cdot
\delta_e,\delta_e)$, where $\{\delta_g\}_{g\in\mb{ F}_n}$ is the
standard orthonormal basis in $l^2(\mb{F}_n)$. Every
$\mc{L}(\mb{F}_n)$ is generated as a von Neumann algebra by a
semicircular system with $n$ elements (\cite{16}). We denote by
$\mc{M}_k^{sa}=\mc{M}_k^{sa}(\mb{C})$ the set of $k\times k$
self-adjoint complex matrices and by $\tau_k$ its unique
normalized trace. $\tau_k$ induces the $2$-norm $||\cdot
||_2:\mc{M}_k^{sa}\rightarrow \mb{R}_+$ and the euclidean norm
$||\cdot ||_e:=\sqrt{k}||\cdot ||_2$. If $B$ is a measurable
subset of a $m$ dimensional (real) manifold then $\mbox{vol}_m(B)$
will denote the Lebesgue measure of $B$. The free entropy $\chi
(x_1,\ldots ,x_n)$ of a finite family of self-adjoint elements was
introduced in \cite{14} but we will recall (\cite{15}) the
definition of the modified free entropy which is better suited for
applications. For self-adjoint elements $x_1,\ldots ,x_{n+m}\in
\mc{M}$ one defines first the set of matricial microstates
\begin{eqnarray}&&\Gamma_R(x_1 ,\ldots ,x_n:x_{n+1},\ldots ,
x_{n+m};p,k,\epsilon ):=\{(A_1,\ldots ,A_n) \in
(\mc{M}_k^{sa})^n|\\\nonumber&& \hspace{.5 cm}\exists
(A_{n+1},\ldots ,A_{n+m})\in (\mc{M}_k^{sa})^m\,\,\mbox{s.t.}\,\,
||A_j||\leq R,\,\, \left|\tau (x_{i_1}\ldots
x_{i_q})\right.\\\nonumber&& \hspace{.5 cm}\left.
-\tau_k(A_{i_1}\ldots A_{i_q})\right| <\epsilon\,\,\forall
j,i_1,\ldots i_q\,\,\in \{1,\ldots ,n+m\}\,\, \forall 1\leq q\leq
p\} \end{eqnarray} where $R,\epsilon >0$ and $p,k\in \mb{ N}$, and
then
\begin{eqnarray}&&\chi_R(x_1 ,\ldots ,x_n:x_{n+1},\ldots
,x_{n+m};p,k,\epsilon )\\\nonumber&&\hspace{.5 cm}= \log
(\mbox{vol}_{nk^2}(\Gamma_R(x_1 ,\ldots ,x_n: x_{n+1},\ldots
,x_{n+m};p,k,\epsilon ))),
\end{eqnarray}
\begin{eqnarray}&&
\chi_R(x_1 ,\ldots ,x_n:x_{n+1},\ldots ,x_{n+m};p,\epsilon
)\\\nonumber &&\hspace{.5 cm}= \limsup_{k\rightarrow\infty
}\left(\frac{1}{k^2}\chi_R(x_1 ,\ldots ,x_n: x_{n+1},\ldots
,x_{n+m};p,k,\epsilon )+\frac{n}{2}\log k\right),
\end{eqnarray}
\begin{eqnarray}&&
\chi_R(x_1 ,\ldots ,x_n:x_{n+1},\ldots
,x_{n+m})\\\nonumber&&\hspace{.5 cm}=\inf\{\chi_R(x_1 ,\ldots
,x_n:x_{n+1},\ldots , x_{n+m};p,\epsilon )| p\in \mb{N}, \epsilon
>0\},
\end{eqnarray}
\begin{eqnarray}&&
\chi (x_1 ,\ldots ,x_n:x_{n+1},\ldots ,x_{n+m})\\\nonumber
&&\hspace{.5 cm}= \sup\{\chi_R(x_1 ,\ldots ,x_n:x_{n+1},\ldots ,
x_{n+m})|R>0\}.
\end{eqnarray}
When taking the last sup it suffices though to assume $0<R\leq
\mbox{max}\{ ||x_1||,$ $\ldots ,||x_{n+m}||\}$ rather than
$0<R<\infty$ (\cite{14}, \cite{15}). The quantity $\chi (x_1 ,$ $
\ldots ,x_n:x_{n+1},\ldots ,x_{n+m})$ is the free entropy of $x_1
,\ldots ,x_n$ in the presence of $x_{n+1},\ldots ,x_{n+m}$. If
$m=0$ it is called the free entropy of $x_1 ,\ldots ,x_n$ and
denoted $\chi (x_1 ,\ldots ,x_n)$. If $\{x_{n+1},\ldots ,$
$x_{n+m}\}\subset \{x_1 ,\ldots$ $,x_n\}''$ then (\cite{15})
$$\chi
(x_1 ,\ldots ,x_n:x_{n+1}, \ldots ,x_{n+m})=\chi (x_1 ,\ldots
,x_n).$$ For a single self-adjoint element $x=x^*\in \mc{M}$ one
has (\cite{14}): $$\chi (x)=\frac{3}{4}+\frac{1}{2}\log 2\pi
+\int\int\log |s-t|d\mu (s)d\mu (t),$$ where $\mu$ is the
distribution of $x$. If $x_1,\ldots ,x_n$ are $n$ self-adjoint
free elements of $\mc{M}$ then $\chi (x_1 ,\ldots ,x_n)=\chi
(x_1)+\ldots +\chi (x_n)$ (\cite{14}). The converse is also true
(\cite{21}), provided that $\chi (x_i)>-\infty\,\,\forall 1\leq
i\leq n$. In particular, the free entropy of a finite semicircular
system is finite, hence the free group factor $\mc{L}(\mb{F}_n)$
has a system of generators with finite free entropy for $2\leq
n<\infty$.
\section{Noncommutative power series and free entropy}\label{chser}
\setcounter{equation}{0} The main result of the present section
states that if a (not necessarily free) family of $m$ self-adjoint
noncommutative random variables with finite free entropy can be
generated as noncommutative power series by another family of $n$
self-adjoint noncommutative random variables, then $n\geq m$. In
other words, a finite system with finite free entropy has minimal
cardinality among all finite systems of self-adjoint elements that
are equivalent under the noncommutative analytic functional
calculus. Thus, we recover D. Voiculescu's result from \cite{14},
with the observation that our approach does not require the
assumption of freeness.

We review first a few facts concerning the theory of systems of
algebraic equations (\cite{20}), necessary in the proof of Lemma
\ref{lema1}. If $g_1,\ldots ,g_n$ are forms in $n$ variables, then
there exists a polynomial (the rezolvent) in their coefficients,
$R(g_1,\ldots ,g_n)$, with the property that $R(g_1,\ldots
,g_n)=0$ if and only if the system $g_1(\xi_1,\ldots
,\xi_n)=\ldots =g_n(\xi_1,\ldots ,\xi_n)=0$ has a nontrivial
solution. If $h_1,\ldots ,h_{n-1}$ are $n-1$ forms in $n$
variables and $h_n(u)(\xi_1,\ldots ,\xi_n):=u_1\xi_1+\ldots
+u_n\xi_n$, then $R_u(h_1,\ldots ,h_{n-1}):=R(h_1,\ldots ,$ $
h_{n-1}, h_n(u))$ (the $u$-rezolvent) is either identically equal
to $0$, or a form of degree $\mbox{deg}(h_1)\cdot\ldots\cdot
\mbox{deg}(h_{n-1})$ in $u=(u_1,\ldots ,u_n)$. In the first case,
the system $h_1=\ldots =h_{n-1}=0$ has infinitely many solutions
$[(\xi_1,\ldots ,\xi_n)]\in \mb{P}\mb{C}^{n-1}$ and in the second,
all the solutions $[(\xi_1,\ldots ,\xi_n)]\in \mb{P}\mb{C}^{n-1}$
are given by the factorization of $R_u(h_1,\ldots ,h_{n-1})$ (and
thus, the system admits at most $\mbox{deg}(h_1)\cdot\ldots\cdot
\mbox{deg}(h_{n-1})$ solutions - B\' ezout's Theorem).

Let $f_1,\ldots ,f_n\in \mb{R}[\Xi_1,\ldots , \Xi_n]$ be $n$
polynomials in $n$ indeterminates, of degrees $d_1,\ldots ,d_n$,
respectively. For $a=(a_1,\ldots ,a_n)\in \mb{R}^n$ define
$$F_{i,a_i}(\xi_1,\ldots ,\xi_{n+1})
=\xi_{n+1}^{d_i}\left(f_i\left(\frac{\xi_1}{\xi_{n+1}},\ldots
,\frac{\xi_n}{\xi_{n+1}} \right)-a_i\right),\forall 1\leq i\leq
n.$$ B\' ezout's Theorem implies that the system of equations
$f_1(\xi_1,\ldots ,\xi_n)$  $=a_1,\ldots , f_n(\xi_1,$ $\ldots
,\xi_n)=a_n$ admits at most $d_1\cdot\ldots\cdot d_n$ solutions
$(\xi_1,\ldots ,\xi_n)\in \mb{C}^n$ if $R_u(F_{1,a_1},\ldots
,F_{n,a_n})\not\equiv 0$. Note also that the set \[S_u(f_1,\ldots
,f_n):=\left\{(a_1,\ldots ,a_n)\in\mb{R}^n\,\,|\,\,
R_u(F_{1,a_1},\ldots ,F_{n,a_n})\not\equiv 0\right\}\] is either
open and dense in $\mb{R}^n$, or empty.

We proceed now with Lemma \ref{lema1} which gives an upper bound
for the Lebesgue measure of the intersection of an algebraically
parameterized manifold embedded in $\mb{R}^m$, with the unit ball
of $\mb{R}^m$. This Lemma will be of further use in estimating the
volumes of various sets of matricial microstates which will appear
as sets of points within given distance from such manifolds.
\begin{lem}\label{lema1}
For integers $n\leq m$ and polynomials $f_1,\ldots ,f_m \in
\mb{R}[\Xi_1,$ $\ldots ,\Xi_n]$ define $f=(f_1,\ldots
,f_m):\mb{R}^n\rightarrow \mb{R}^m$. If the polynomials
$\det\left(\frac{\partial f_J}{\partial\xi}\right)$ are not
identically equal to $0$ $\forall J\in \lbrace (i_1,\ldots
,i_n)|1\leq i_1<\ldots <i_n\leq m\rbrace$ and if
$S_u=S_u(f_1,\ldots ,f_n)\not=\emptyset$, then
\begin{equation}\label{eq1}
\int_{f^{-1}\left(\overline{B(0,1)}\right)} \left( \sum_{|J|=n}
{\det}^2\left(\frac{\partial f_J}{\partial\xi}\right)
\right)^{\frac{1}{2}} d\xi \leq\binom{m}{n}\cdot C \cdot
\mbox{vol}_n (B(0,1))\,\,,
\end{equation}
where $C=C(\deg (f))=\max\{\deg (f_{i_1})\cdot\ldots\cdot \deg
(f_{i_n})\,\,|\,\,1\leq i_1<\ldots <i_n\leq m\}$ and
$B(0,1)=B_n(0,1)$ is the unit ball in $\mb{R}^n$.
\end{lem}
\begin{proof} We consider first the case $m=n$. Let $S$ denote the set of all
irregular values of $f$, $S=f\left( \lbrace \xi\in
\mb{R}^n\,\,|\,\, \mbox{rank}(df_\xi) < n\rbrace \right)$. It
suffices to show that (\ref{eq1}) holds with
$f^{-1}\left(\overline{B(0,1)} \setminus S_{\epsilon}\right)$
replacing $f^{-1}\left(\overline{B(0,1)}\right)$, where
$S_{\epsilon}$ is an arbitrary open set that contains $S\cup
(\mb{R}^n\setminus S_u)$. For any $a=(a_1,\ldots , a_n)\in
\mbox{Range}(f)\cap \overline{B(0,1)}\setminus S_\epsilon $ the
set $f^{-1}(\{a\})$ has at most $C=\mbox{deg}(f_1)\cdot\ldots\cdot
\mbox{deg}(f_n)$ elements, say $f^{-1}(\{ a\} )=\{b_1,\ldots
,b_{p(a)}\}$ for some $1\leq p(a)\leq C$. There exist an open ball
$B_a\ni a$ and open neighborhoods $V_1^a\ni b_1,\ldots
,V_{p(a)}^a\ni b_{p(a)}$ such that $B_a$ and $V_i^a$ are
diffeomorphic via $f$ for $1\leq i\leq p(a)$ and
$f^{-1}(B_a)=\cup_{i=1}^{p(a)}V_i^a$. Since it is compact, we can
cover $\mbox{Range}(f)\cap\overline{B(0,1)} \setminus
S_{\epsilon}$ with a finite set of such open balls $B_{a_1},\ldots
,B_{a_k}$. This covering determines a finite partition of
$\mbox{Range}(f)\cap\overline{B(0,1)} \setminus S_{\epsilon}$, say
$W_1,\ldots ,W_t$. For each $1\leq j\leq t$ choose a unique $1\leq
l=l(j)\leq k$ such that $W_j\subset B_{a_l}$ and
$f^{-1}(W_j)=T_{j1}\cup\ldots \cup T_{jp(a_l)}$ where
$T_{ji}\subset V_i^{a_l}$ and $W_j$ and $T_{ji}$ are diffeomorphic
via $f$ for all $1\leq i\leq p(a_l)$.
\begin{eqnarray}&&\int_{f^{-1}\left(\overline{B(0,1)} \setminus
S_{\epsilon}\right)} \left|\mbox{det}\left(\frac{\partial
f}{\partial \xi}\right)\right| d\xi =
\sum_{j=1}^{t}\int_{f^{-1}(W_j)}\left|\mbox{det}\left(\frac{\partial
f}{\partial \xi}\right)\right| d\xi\\\nonumber&&\hspace{1 cm}
=\sum_{j=1}^{t}\sum_{i=1}^{p(a_{l(j)})}
\int_{T_{ji}}\left|\mbox{det}\left(\frac{\partial f}{\partial
\xi}\right)\right|d\xi
=\sum_{j=1}^{t}\sum_{i=1}^{p(a_{l(j)})}\mbox{vol}_n (W_j)
\\\nonumber &&\hspace{1 cm}\leq C\sum_{j=1}^{t}\mbox{vol}_n(W_j)=C\cdot \mbox{vol}_n
\left(\overline{B(0,1)}\setminus S_{\epsilon}\right).
\end{eqnarray}
In the case $m>n$ one has the following estimates:
\begin{eqnarray}\nonumber&&\hspace{-.5 cm}\int_{f^{-1}\left(\overline{B(0,1)}\right)} \left(
\sum_{|J|=n} \mbox{det}^2 \left(\frac{\partial f_J}{\partial
\xi}\right) \right)^{\frac{1}{2}} d\xi
\leq\int_{f^{-1}\left(\overline{B(0,1)} \right)}
\sum_{|J|=n}\left| \mbox{det}\left(\frac{\partial f_{J}}{\partial
\xi}\right)\right| d\xi\\\nonumber&&\hspace{.5 cm}\leq
\sum_{|J|=n}\int_{f_J^{-1}\left(\overline{B(0,1)} \right)} \left|
\mbox{det}\left(\frac{\partial f_{J}}{\partial
\xi}\right)\right|d\xi \leq \binom{m}{n} \cdot C \cdot
\mbox{vol}_n (B(0,1))\,.
\end{eqnarray}
\end{proof}
Lemma \ref{lema1} will be used in the proof of Proposition
\ref{prop1}. The $k\times k$ matricial microstates of $x_1,\ldots
,x_m$ are points within euclidean distance $2\omega\sqrt{mk}$ from
the range of a polynomial function in the matricial microstates of
$y_1,\ldots ,y_n$ provided that each $x_i$ is within $||\cdot
||_2$-distance $\omega$ from noncommutative polynomials in
$y_1,\ldots ,y_n$.
\begin{prop}\label{prop1}
Let $P_1,\ldots ,P_m\in \mb{C}<Y_1,\ldots,Y_n>$ be complex
polynomials in $n$ noncommutative self-adjoint variables. Assume
that $(\mc{M},\tau )$ is a $\mbox{\!I\!I}_1$-factor and
$\{x_1,\ldots,x_m\}\subset \mc{M}$ is a finite set of self-adjoint
generators of $\mc{M}$. If $\{y_1,\ldots,y_n\}\subset \mc{M}$ is
another finite set of self-adjoint generators of $\mc{M}$ with
$n<m$ and such that
$$ ||x_i-P_i(y_1,\ldots,y_n)||_2<\omega\,\,\forall \,\,1\leq i\leq m$$
for some positive constant $\omega\in (0,a]$, then
\begin{equation}
\chi (x_1,\ldots,x_m)\leq C(m,n,a)+(m-n)\log \omega +n\log d
\end{equation}
where $C(m,n,a)$ is a constant that depends only on $a=\max\lbrace
||x_1||_2+1,\ldots, ||x_m||_2+1\rbrace$, $m$, $n$ and
$d=\max\lbrace \deg (P_1),\ldots , \deg (P_m)\rbrace$.
\end{prop}
\begin{proof} Eventually replacing each $P_i$ by $\frac{1}{2}(P_i+P_i^*)$ we can
assume from the beginning that $P_i=P_i^*\,\,\forall 1\leq i\leq
m$. For $R>0$, integer $p\geq 1$ and $\epsilon >0$ consider
$$(A_1,\ldots,A_m)\in \Gamma_R(x_1,\ldots,x_m:y_1,\ldots,y_n;p,k,\epsilon)\,\,.
$$
If $p$ is large enough and $\epsilon >0$ is sufficiently small,
then one can find matrices $B_1,\ldots,B_n\in \mc{M}_k^{sa}$ such
that $||B_1||,\ldots,||B_n||\leq R$ and
$$||A_i-P_i(B_1,\ldots,B_n)||_2<\omega\,\,\forall 1\leq i\leq m$$
or equivalently,
$$||A_i-P_i(B_1,\ldots,B_n)||_e<\omega\sqrt{k}\,\,\forall 1\leq i\leq m\,\,.
$$

With the identifications $g=(g_1,\ldots
,g_{mk^2}):(\mc{M}_k^{sa})^n\cong \mb{R}^{nk^2}$ $\rightarrow
(\mc{M}_k^{sa})^m\cong \mb{R}^{mk^2}$, $(B_1,\ldots,B_n)=
(\xi_1,\ldots ,\xi_{nk^2})\in \mb{R}^{nk^2}$, $g(B_1,\ldots,B_n)=
(P_1(B_1,\ldots,$ $B_n),\ldots ,P_m(B_1,\ldots,B_n))$, the
previous inequality becomes
$$\left|\left|(A_i)_{1\leq i\leq m}-g(\xi_1,\ldots ,\xi_{nk^2})\right|\right|_e
<\omega\sqrt{mk}.$$ At the cost of introducing an additional
variable $\xi_{nk^2+1}\in\mb{R}$, we can assume that the
components of $g$ are $mk^2$ homogeneous polynomial functions in
the variables $\xi_1,\ldots ,\xi_{nk^2+1}$, of degrees $\leq d$.

Let now $f_1,\ldots ,f_{mk^2}$ be arbitrary homogeneous polynomial
functions in $\xi_1,\ldots ,\xi_{nk^2+1}$ such that
$\mbox{deg}(f_j)=\mbox{deg}(g_j)\,\, \forall 1\leq j\leq mk^2$.
For every multiindex $J=(j_1,\ldots ,j_{nk^2+1})$ with $1\leq
j_1<\ldots <j_{nk^2+1}\leq mk^2$, $S_u(f_{j_1},\ldots
,f_{j_{nk^2+1}})=\emptyset$ is equivalent to the fact that the
coefficients of $f_{j_1},\ldots ,f_{j_{nk^2+1}}$ satisfy a certain
system of algebraic equations. Hence the set
\begin{eqnarray}&&\Omega_1=\{f=(f_1,\ldots
,f_{mk^2})\,\,|\,\,\mbox{deg}(f_j)=\mbox{deg}(g_j)\,\, \forall
1\leq j\leq mk^2,\\\nonumber&& \hspace{1 cm}S_u(f_{j_1},\ldots
,f_{j_{nk^2+1}}) \not=\emptyset\,\,\forall J=(j_1,\ldots
,j_{nk^2+1})\}
\end{eqnarray} is open and dense in its natural ambient linear
space. Similarly, the set
\begin{eqnarray}&&\Omega_2=\bigg\{f=(f_1,\ldots
,f_{mk^2})\,\,|\,\,\mbox{deg}(f_j)=\mbox{deg}(g_j)\,\, \forall
1\leq j\leq mk^2,\\\nonumber&& \hspace{1 cm}
\mbox{det}\left(\frac{\partial f_J}{\partial \xi}\right)
\not\equiv 0\,\,\forall J=(j_1,\ldots ,j_{nk^2+1})\bigg\}
\end{eqnarray}
is also open and dense in the same linear space..

The matrix $df_\xi$ has $\binom{mk^2}{nk^2+1}$ minors of dimension
$(nk^2+1)\times (nk^2+1)$ and all these minors have a nontrivial
common zero only if (\cite{20}) a certain system of algebraic
equations in the coefficients of $f_1,\ldots ,f_{mk^2}$ has a
solution. Moreover, not all the polynomials appearing in this
system are identically equal to $0$. It follows that the set
\begin{eqnarray}&&\Omega_3=\big\{f=(f_1,\ldots
,f_{mk^2})\,\,|\,\,\mbox{deg}(f_j)=\mbox{deg}(g_j)\,\,\forall
1\leq j\leq mk^2,\\\nonumber&&\hspace{1 cm}
\mbox{rank}(df_\xi)=nk^2+1\,\,\forall \xi\in
\mb{R}^{nk^2+1}\setminus\{0\} \big\} \end{eqnarray} contains a
subset which is open and dense in the linear space previously
considered. Therefore there exists an element $f\in
\Omega_1\cap\Omega_2\cap\Omega_3$ such that
$\left|\left|f(\xi_1,\ldots ,\xi_{nk^2+1})-g(\xi_1,\ldots
,\xi_{nk^2+1})\right|\right|_e <\omega\sqrt{mk}$ $\forall
|\xi_i|\leq R$ $\forall 1\leq i\leq nk^2+1$, hence
$\left|\left|(A_i)_{1\leq i\leq m}-f(\xi_1,\ldots
,\xi_{nk^2+1})\right|\right|_e <2\omega\sqrt{mk}$. The function
$f$ satisfies the hypothesis of Lemma \ref{lema1} and its
components are homogeneous polynomials. Moreover, it has the
property that $\mbox{dist}_e \left((A_i)_{1\leq i\leq
m},\,\,\mbox{Range}\left(f\right)\right) <2\omega\sqrt{mk}$ and it
does not depend on the system $(A_i)_{1\leq i\leq m}$.

We have $||(A_1,\ldots,A_m)||_e\leq
 a\sqrt{mk}$ (if $\epsilon >0$ is small enough) hence the set
 of matricial microstates $(A_1,\ldots,A_m)$ of $(x_1,\ldots,x_m)$ such that
$\mbox{dist}_e((A_1,$ $\ldots,A_m), \mbox{Range}(f))<2\omega
\sqrt{mk}$ is contained in the $(mk^2,nk^2$ $+1)$-tube of radius
$2\omega \sqrt{mk}$ around $\mbox{Range}(f)\cap
B_{mk^2}(0,(a+2\omega)\sqrt{mk})$. If $B$ is a small ball in
$\mb{R}^{nk^2+1}\setminus\{0\}$ and if $V_B(2\omega\sqrt{mk})$
denotes the $(mk^2,nk^2$ $+1)$-tube of radius $2\omega \sqrt{mk}$
around $f(B)$, then the formula for volumes of tubes (\cite{22})
implies
\begin{eqnarray}&&
\mbox{vol}_{mk^2}(V_B(2\omega\sqrt{mk}))=\mbox{vol}_{mk^2-nk^2-1}
(B_{mk^2-nk^2-1}(0,1))\\\nonumber&& \cdot\sum_{\overset{e\equiv
0\,\,mod\,\,2}{0\leq e\leq nk^2+1}}\frac{(2\omega\sqrt{mk})^{e+
mk^2-nk^2-1}k_{B,e}}{(mk^2-nk^2+1)(mk^2-nk^2+3)\ldots
(mk^2-nk^2-1+e)}. \end{eqnarray} With the notations from \cite{22}
one has $k_{B,e}=\int_{f(B)}H_eds$ and
$$H_e=\frac{1}{2^e(e/2)!}\sum_{\sigma\in\Sigma_e}\mbox{sgn} (\sigma
) \sum_{\alpha_1,\ldots
,\alpha_e=1}^{nk^2+1}H_{\alpha_1\alpha_2}^{ \alpha_{\sigma
(1)}\alpha_{\sigma (2)}}H_{\alpha_3\alpha_4}^{\alpha_{ \sigma
(3)}\alpha_{\sigma (4)}}\ldots$$ where
$H_{\alpha\beta}^{\lambda\mu}$ denotes the Riemann tensor of
$f(B)$. Assuming without loss of generality that $\mbox{deg}(f_j)=
d$ $\forall 1\leq j\leq mk^2$, one can verify that each
$H_{\alpha\beta}^{\lambda\mu}(f(\xi ))$ is a sum of quotients of
homogeneous polynomials where all numerators have degree
$6(d-1)(nk^2+1)-2d$ and all denominators have degree
$6(d-1)(nk^2+1)$, hence $H_e$ is a rational function in $\xi$ and
in the coefficients of $f(\xi)$. Due to its intrinsic nature,
$H_e$ is independent of the embedding of $\mbox{Range}(f)$ in
$\mb{R}^{mk^2+1}$, in particular it is invariant under orthogonal
transformations in $\mb{R}^{mk^2+1}$. Since there exist
sufficiently many polynomials $f(\xi)$ such that $\mbox{Range}(f)$
is flat, this entails $H_e=0$ $\forall 2\leq e\leq nk^2+1$,
$e\equiv 0\,\,\mbox{mod}\,\,2$. Therefore the volume of the
$(mk^2,nk^2$ $+1)$-tube of radius $2\omega \sqrt{mk}$ around
$f(B)$ is
$\mbox{vol}_{mk^2}(V_B(2\omega\sqrt{mk}))=\mbox{vol}_{mk^2-nk^2-1}
(B_{mk^2-nk^2-1}(0,1))\cdot (2\omega\sqrt{mk})^{
mk^2-nk^2-1}\cdot\int_{f(B)}ds$ and with Lemma \ref{lema1} and
inequality \begin{equation}\label{g}\frac{1}{\Gamma
\left(1+\frac{nk^2+1}{2}\right)}\cdot \frac{1}{\Gamma
\left(1+\frac{mk^2-nk^2-1}{2}\right)}\leq\frac{2^{\frac{mk^2}{2}}}
{\Gamma\left(1+\frac{mk^2}{2}\right)}\end{equation} we obtain the
following estimate:
\begin{eqnarray}&&\mbox{vol}_{mk^2}(\Gamma_R(x_1,\ldots,x_m:y_1,\ldots,
y_n;p,k,\epsilon))\leq \binom{mk^2}{nk^2+1}\cdot
C(d)\\\nonumber&&\hspace{.5 cm}\cdot
\mbox{vol}_{nk^2+1}\left(B(0,(a+2\omega)\sqrt{mk})\right)\cdot
\mbox{vol}_{mk^2-nk^2-1}(B(0,1))\\\nonumber&&\hspace{.5 cm}\cdot
(2\omega \sqrt{mk})^{mk^2-nk^2-1}= \binom{mk^2}{nk^2+1}\cdot
C(d)\cdot\pi^{\frac{nk^2+1}{2}}
\\\nonumber&&\hspace{.5 cm}\cdot
\frac{(a+2\omega)^{nk^2+1}(mk)^{\frac{nk^2+1}{2}}
\pi^{\frac{mk^2-nk^2-1}{2}}(2\omega)^{mk^2-nk^2-1}(mk)^
{\frac{mk^2-nk^2-1}{2}}}
{\Gamma\left(1+\frac{nk^2+1}{2}\right)\Gamma\left(1+\frac{mk^2-nk^2-1}{2}
\right)}\\\nonumber&&\hspace{.5 cm}\leq \binom{mk^2}{nk^2+1}\cdot
C(d)\cdot
\frac{\pi^{\frac{mk^2}{2}}(mk)^{\frac{mk^2}{2}}2^{\frac{mk^2}{2}}
(3a)^{nk^2+1}(2\omega
)^{mk^2-nk^2-1}}{\Gamma\left(1+\frac{mk^2}{2} \right)}.
\end{eqnarray}
The last inequality implies further
\begin{eqnarray}&&\chi_R(x_1,\ldots,x_m:y_1,\ldots,y_n;p,k,\epsilon)\leq
\frac{1}{k^2}\log \binom{mk^2}{nk^2+1} +\frac{1}{k^2}\log
C(d)\\\nonumber&&\hspace{.5 cm} +\frac{m}{2}\log \pi
+\left(\frac{3m}{2}-n\right)\log 2 +n\log (3a)+ \frac{m}{2}\log
(mk)\\\nonumber&&\hspace{.5 cm}+(m-n)\log\omega
-\frac{1}{k^2}\log\Gamma\left(1+\frac{mk^2}{2}\right)+\frac{m}{2}\log
k+o(1).
\end{eqnarray}
Note that one has
$\frac{1}{k^2}\log\Gamma\left(1+\frac{mk^2}{2}\right)=\frac{m}{2}\log
\frac{mk^2}{2e}+o(1)$, $C(d)\leq d^{nk^2+1}$ and
$\frac{1}{k^2}\log \binom{mk^2}{nk^2+1}= m\log m-n\log n-(m-n)\log
(m-n)+o(1)$, therefore
\begin{eqnarray}&&\chi_R(x_1,\ldots,x_m:y_1,\ldots,y_n;p,k,\epsilon)\leq m\log
m-n\log n+n\log d\\\nonumber&&\hspace{.5 cm}-(m-n)\log
(m-n)+\frac{m}{2}\log\pi +\left(\frac{3m}{2}-n\right)\log 2 +n\log
(3a)\\\nonumber&&\hspace{.5 cm}+\frac{m}{2}\log m+\frac{m}{2}\log
k+(m-n)\log\omega -\frac{m}{2} \log \frac{m}{2e} -m\log
k\\\nonumber&&\hspace{.5 cm}+\frac{m}{2}\log k
+o(1)=C(m,n,a)+(m-n)\log\omega +n\log d+o(1). \end{eqnarray} By
taking the appropriate limits after $k,p,\epsilon$, we finally
obtain
$$\chi_R(x_1,\ldots,x_m:y_1,\ldots,y_n)\leq
C(m,n,a)+(m-n)\log\omega +n\log d ,$$ and since $R>0$ is
arbitrary, $\chi(x_1,\ldots,x_m:y_1,\ldots,y_n)\leq
C(m,n,a)+(m-n)\log\omega +n\log d$. Recall now that
$\{x_1,\ldots,x_m\}$ is a system of generators of $\mc{M}$, hence
$\chi (x_1,\ldots,x_m)=\chi(x_1,\ldots,x_m:y_1,\ldots,y_n)$.
\end{proof}
Let $Y_1,\ldots,Y_n$ be noncommutative indeterminates and let
$$P(Y_1,\ldots,Y_n)=\sum_{k=0}^{\infty}\sum_{1\leq i_1,\ldots ,i_k\leq n}
a_{i_1\ldots i_k} Y_{i_1}\ldots Y_{i_k}$$ be a noncommutative
power series in $Y_1,\ldots,Y_n$, with complex coefficients.
Following \cite{14}, we say that $R>0$ is a radius of convergence
of $P$ if
$$\sum_{k=0}^{\infty}\sum_{1\leq i_1,\ldots ,i_k\leq n}
|a_{i_1\ldots i_k}| R^k<\infty\,\,.$$ It is well-known from the
theory of power series that if $0<R_0<R$, then
$$\sum_{k=q+1}^{\infty}\sum_{1\leq i_1,\ldots ,i_k\leq n}
|a_{i_1\ldots i_k}|
R_0^k=O\bigg(\bigg(\frac{R_0}{R}\bigg)^{q+1}\bigg)\,\,.$$ Theorem
\ref{prop2} is basically Corollary 6.12 in \cite{14}, with the
observation that the freeness of $\{x_1,\ldots ,x_m\}$ assumed
there has been dropped.
\begin{thm}\label{prop2}
Let $x_1,\ldots ,x_m$ and $y_1,\ldots ,y_n$ be self-adjoint
noncommutative random variables in a $\mbox{\!I\!I}_1$-factor
$(\mc{M},\tau )$ such that $y_1,\ldots ,y_n\in\{x_1,\ldots ,$
$x_m\}''$ and $\chi (x_1,\ldots,x_m)>-\infty$. If
$x_i=P_i(y_1,\ldots ,y_n)\,\,\forall 1\leq i\leq m$, where
$(P_i)_{1\leq i\leq m}$ are noncommutative power series having a
common radius of convergence $R> b=\max\{||y_1||,\ldots
,||y_n||\}$, then $n\geq m$.
\end{thm}
\begin{proof} Suppose that $m>n$. For $1\leq i\leq m$, $x_i$ is a
noncommutative power series of $y_1,\ldots ,y_n$ i.e.,
$$x_i=\sum_{k=0}^{\infty}\sum_{1\leq i_1,\ldots ,i_k\leq n}
a_{i_1\ldots i_k}^{(i)} y_{i_1}\ldots y_{i_k}\,\,.$$ For every
integer $q\geq 0$, $P_{i,q}(y_1,\ldots ,y_n):=\sum_{k=0}^q
\sum_{1\leq i_1, \ldots ,i_k\leq n} a_{i_1\ldots i_k}^{(i)}
y_{i_1}$ $\ldots y_{i_k}$ is a noncommutative polynomial of degree
$\leq q$ and moreover
\begin{eqnarray}&&||x_i-P_{i,q}(y_1,\ldots ,y_n)||_2=
\bigg|\bigg|\sum_{k=q+1}^{\infty} \sum_{1\leq i_1, \ldots ,i_k\leq
n} a_{i_1\ldots i_k}^{(i)} y_{i_1}\ldots
y_{i_k}\bigg|\bigg|_2\\\nonumber &&\hspace{1 cm}\leq
\sum_{k=q+1}^{\infty} \sum_{1\leq i_1, \ldots ,i_k\leq n}
|a_{i_1\ldots i_k}^{(i)}|
b^k=O\bigg(\bigg(\frac{b}{R}\bigg)^{q+1}\bigg).
\end{eqnarray} The estimate of free entropy
from Proposition \ref{prop1} implies $\chi (x_1,\ldots,x_m)$ $\leq
C(m,n,a)+(m-n)\log\left(\frac{b}{R}\right)^{q+1} +n\log q +O(1)$
and letting $q$ tend to $\infty$, one obtains that $\chi
(x_1,\ldots,x_m)=-\infty$, contradiction. \end{proof}

Let $\mc{N}$ be a $*$-algebra in a $W^*$-probability space
$(\mc{M},\tau)$. Suppose that $\mc{N}$ is finitely generated and
let $\{x_1,\ldots ,x_m\}$ be a system of self-adjoint generators.
Let also $\{y_1,\ldots ,y_n\}$ be another set of self-adjoint
elements that generate $\mc{N}$ {\it algebraically} as a
$*$-algebra. In particular, there exist noncommutative polynomials
$(P_i)_{1\leq i\leq m}$ such that $x_i=P_i(y_1,\ldots
,y_n)\,\,\forall 1\leq i\leq m$. In this context, Corollary
\ref{cor1} is an immediate consequence of Theorem \ref{prop2}.
\begin{cor}\label{cor1}
If $\chi (x_1,\ldots,x_m)>-\infty$ and $*$-alg$\{y_1,\ldots
,y_n\}= *$-alg$\{x_1,$ $\ldots ,x_m\}$ then $n\geq m$, so any $2$
systems of self-adjoint elements with finite free entropy that
generate $\mc{N}$ algebraically as a $*$-algebra have the same
cardinality.
\end{cor}
D. Voiculescu proved in \cite{17} that the modified free entropy
dimension (\cite{15}) of a finite set of self-adjoint elements
that generate algebraically a $*$-algebra $\mc{N}$ is independent
of the set of generators. Recall (\cite{15}) the definition of the
modified free entropy dimension: $$\delta_0(x_1,\ldots ,
x_m)=m+\limsup_{\omega\rightarrow 0} \frac{\chi (x_1+\omega
s_1,\ldots ,x_1+\omega s_m:s_1,\ldots ,s_m)} {|\log\omega |},$$
where $\{s_1,\ldots ,s_m\}$ is a semicircular system free from
$\{x_1,\dots ,x_m\}$. One has $\delta_0(x_1,\ldots ,x_m)\leq m$ in
general, and also $0\leq \delta_0(x_1,\ldots ,$ $x_m)$ if
$\{x_1,\ldots ,x_m\}$ $\subset \mc{L}(\mb{F}_p)$ for some $p$.
Considering two sets $\{x_1,\ldots ,x_m\}$ and $\{y_1,\ldots
,y_n\}$ of self-adjoint elements that generate algebraically the
$*$-algebra $\mc{N}$ and noticing that $\{y_1,\ldots
,y_n\}\subset\{x_1+\omega s_1,\ldots ,x_1+\omega s_m,s_1,\ldots
,s_m\}''$, one has \begin{eqnarray}&&\delta_0(x_1,\ldots
,x_m)=m\\\nonumber&& \hspace{.5 cm}+\limsup_{\omega\rightarrow 0}
\frac{\chi (x_1+\omega s_1,\ldots ,x_1+\omega s_m:s_1,\ldots
,s_m,y_1, \ldots ,y_n)} {|\log\omega |}\\\nonumber&& \hspace{.5
cm} \leq m+\limsup_{\omega\rightarrow 0} \frac{\chi (x_1+\omega
s_1,\ldots ,x_1+\omega s_m:y_1, \ldots ,y_n)} {|\log\omega |}.
\end{eqnarray}
Also, $||x_i+\omega s_i-P_i(y_1,\ldots ,y_n)||=||\omega
s_i||\leq\omega\,\, \forall 1\leq i\leq m$, and with Proposition
\ref{prop1} we obtain \begin{eqnarray}&&\delta_0(x_1,\ldots
,x_m)\leq m\\\nonumber&&\hspace{.5 cm}+\limsup_{\omega\rightarrow
0} \frac{C(m,n,a)+(m-n)\log\omega +n\log d} {|\log\omega |}\leq
m+n-m=n, \end{eqnarray} where $a=\max\lbrace ||x_1||_2+1,
\ldots,||x_m||_2+1,||y_1||_2+1, \ldots,||y_n||_2+1\rbrace$ and
$d=\max\{\mbox{deg}(P_i)\,|\,1\leq i\leq m\}$. In particular, if
there exists a set $\{y_1,\ldots ,y_n\}$ with $\delta_0(y_1,\ldots
,y_n)=n$ which generates $\mc{N}$ algebraically, then $\sup
\{\delta_0(x_1,\ldots ,$ $x_m)\,\,|\,\,*$-alg$\{x_1,\ldots ,x_m\}=
\mc{N}\}=n$.
\section{Indecomposability over nonprime subfactors}\label{chhyp}
\setcounter{equation}{0} In this section we prove that the free
group factor $\mc{L}(\mb{F}_n)$ does not admit an asymptotic
decomposition of the form
$$\lim_{\omega\rightarrow 0}\,^{||\cdot ||_2}\sum_{\overset{1\leq j_1,\ldots
,j_{t+1}\leq f}{1\leq t\leq d}} \mc{N}_{j_1}^\omega \mc{Z}^\omega
\mc{N}_{j_2}^\omega \mc{Z}^\omega\ldots \mc{N}_{j_t}^\omega
\mc{Z}^\omega \mc{N}_{j_{t+1}}^\omega ,$$ where
$\{\mc{Z}^\omega\subset \mc{L}(\mb{F}_n)\}_\omega$ are subsets
with $p$ self-adjoint elements,  $\{\mc{N}_1^\omega ,\ldots ,$ $
\mc{N}_f^\omega\}$ are nonprime subfactors of $\mc{L}(\mb{F}_n)$,
$d\geq 1$ is an arbitrary integer, and $n \geq p+2f +1$. A
nonprime $\mbox{\!I\!I}_1$-factor is just a factor isomorphic to
the tensor product of two factors of type $\mbox{\!I\!I}_1$. For
free group subfactors one has the following: if $n\geq p+2f+2$ and
$\mc{P}\subset \mc{L}(\mb{F}_n)$ is a subfactor of finite index,
then $\mc{P}$ does not admit such an asymptotic decomposition
either. In particular, the hyperfinite dimension of
$\mc{L}(\mb{F}_n)$ is $\geq [\frac{n-2}{2}]+1$ and the hyperfinite
dimension of $\mc{P}$ is $\geq [\frac{n-3}{2}]+1$. For $n=\infty$
this settles a conjecture of L. Ge and S. Popa (\cite{5}): the
hyperfinite dimension of free group factors is infinite. The
definitions of hyperfinite dimension and of asymptotic
decomposition over nonprime subfactors are given next.
\begin{defn}\label{hd}
(\cite{5}) If $\mc{M}$ is a type $\mbox{\!I\!I}_1$-factor, then
the hyperfinite dimension of $\mc{M}$, denoted $\ell_h(\mc{M})$,
is by definition the smallest positive integer $f\in\mb{N}$ with
the property that there exist hyperfinite subalgebras
$\mc{R}_1,\ldots ,\mc{R}_f\subset \mc{M}$ such that
$\overline{\s}^w\mc{R}_1\mc{R}_2\ldots \mc{R}_f=\mc{M}$. If there
is no such positive integer $f$, then by definition,
$\ell_h(\mc{M})=+\infty$.
\end{defn}
\begin{defn}\label{hde}
A type $\mbox{\!I\!I}_1$-factor $\mc{M}$ admits an asymptotic
decomposition over nonprime subfactors, denoted
$$\lim_{\omega\rightarrow 0}\,^{||\cdot ||_2}\sum_{\overset{1\leq j_1,\ldots
,j_{t+1}\leq f}{1\leq t\leq d}} \mc{N}_{j_1}^\omega \mc{Z}^\omega
\mc{N}_{j_2}^\omega \mc{Z}^\omega\ldots \mc{N}_{j_t}^\omega
\mc{Z}^\omega \mc{N}_{j_{t+1}}^\omega\,,$$ provided that $\forall
n\geq 1$ $\forall x_1,\ldots ,x_n\in \mc{M}$ $\forall\omega >0$
$\exists \mc{N}_1^\omega = \mc{N}_1(x_1,\ldots ,x_n;$
$\omega),\ldots ,\mc{N}_f^\omega =\mc{N}_f(x_1,\ldots
,x_n;\omega)$ nonprime subfactors of $\mc{M}$ $\exists
\mc{Z}^\omega =\mc{Z}(x_1,$ $\ldots ,x_n;\omega)\subset \mc{M}$
containing $p$ self-adjoint elements, such that
$$\dist_{||\cdot ||_2}\left(x_j,\sum_{\overset{1\leq j_1,\ldots ,j_{t+1}\leq f}{1\leq
t\leq d}}
\mc{N}_{j_1}^\omega\mc{Z}^\omega\mc{N}_{j_2}^\omega\mc{Z}^\omega\ldots
\mc{N}_{j_t}^\omega\mc{Z}^\omega\mc{N}_{j_{t+1}}^\omega\right)<\omega\,\,
\forall 1\leq j\leq n.$$
\end{defn}
If $\mc{L}(\mb{F}_n)$ admitted an asymptotic decomposition over
nonprime subfactors (Definition \ref{hde}), then the situation
described in Proposition \ref{prop3} (with $\mc{M}=\mc{L}(\mb{
F}_n)$) would take place for an arbitrary $\omega >0$, since any
$\mbox{\!I\!I}_1$-factor is generated by its projections of given
trace ($\frac{1}{2}$, for example):
\begin{lem}
(\cite{50}) Any type $\mbox{\!I\!I}_1$-factor $\mc{M}$ with
separable predual is generated by a countable family of
projections of given trace.
\end{lem}
\begin{proof} Every $\mbox{\!I\!I}_1$-factor with separable predual is generated
by a countable family of abelian subalgebras, so there exist
$\mc{A}_1,\mc{A}_2,\ldots$ abelian subalgebras of $\mc{M}$ that
generate $\mc{M}$ as a von Neumann algebra. If necessary, one can
replace each $\mc{A}_n$ by a maximal abelian subalgebra of
$\mc{M}$ which contains it, hence we can assume that $\mc{A}_n$ is
a maximal abelian subalgebra of $\mc{M}$ $\forall 1\leq n<\infty$.
Being a maximal abelian subalgebra of a type
$\mbox{\!I\!I}_1$-factor, $\mc{A}_n$ has no atoms and thus it is
generated by a countable subset of projections of given trace,
$\forall 1\leq n<\infty$.
\end{proof}
\begin{prop}\label{prop3}
Let $z_1,\ldots ,z_p$ be self-adjoint elements of a
$\mbox{\!I\!I}_1$-factor $\mc{M}$ and let $(\mc{N}_v)_{1\leq v\leq
f}$ be a family of subfactors of $\mc{M}$. Assume that
$\mc{N}_v=\mc{R}^{(v)}_1\vee
\mc{R}^{(v)}_2\simeq\mc{R}^{(v)}_1\otimes \mc{R}^{(v)}_2$ where
$\mc{R}^{(v)}_1,\mc{R}^{(v)}_2$ are $\mbox{\!I\!I}_1$-factors and
assume that $x_1,\ldots ,x_n$ are self-adjoint generators of
$\mc{M}$. Assume moreover that there exist projections of trace
$\frac{1}{2}$, $p^{(v)}_1,\ldots,p^{(v)}_{r_v}\in \mc{R}^{(v)}_2$,
$q^{(v)}_1,\ldots,q^{(v)}_{s_v}\in \mc{R}^{(v)}_1$ and complex
 noncommutative polynomials $(\phi_j)_{1\leq j\leq n}$
 of degrees $\leq d$ (where $d\geq 1$ is fixed)
 in the variables $(z_u)_{1\leq u\leq p}$
 such that
 \begin{equation}\label{cond1}
 \left |\left |x_j-\phi_j\left((p^{(v)}_i)_{\overset{1\leq i\leq r_v}{ 1\leq v\leq f}},
 (q^{(v)}_l)_{\overset{1
 \leq l\leq s_v}{1\leq v\leq f}},
(z_u)_{1\leq
 u\leq p}\right) \right |\right |_2<\omega ,\,\,j=1,\ldots,n,
 \end{equation}
 where $\omega\in (0,a]$ is a given positive number, and such that in all the monomials
 of each
 $\phi_j$ the projections $p^{(v)}_i,q^{(v)}_l$ and $p^{(w)}_k,q^{(w)}_s$ are separated
 by some $z_u$ if $v\ne w$.
 Then
 \begin{equation}\label{estimate1}
\chi(x_1,\ldots,x_n)\leq C(n,p,a,d,f)+(n-p-2f)\log\omega ,
\end{equation}
where $a=\max\lbrace ||x_j||_2+1 |1\leq j\leq n\rbrace$ and
$C(n,p,a,d,f)$ is a constant that depends only on $n,p,a,d,f$.
\end{prop}
\begin{proof} All variables involved are self-adjoint so we can assume that
$\phi_j=\phi_j^*$ $\forall 1\leq j\leq n$. Fix an integer $k_0\geq
1$ and let $R>0$. Let $\mc{M}_{k_0}(\mb{C})\cong \mc{M}^{(v)}_1
\subset \mc{R}^{(v)}_1$, $\mc{M}_{k_0}(\mb{C})\cong
\mc{M}^{(v)}_2\subset \mc{R}^{(v)}_2$ and $\lbrace e^{(v)}_{jl}
\rbrace_{j,l}$ , $\lbrace f^{(v)}_{jl}\rbrace_{j,l}$ be matrix
units for $\mc{M}^{(v)}_1$ and $\mc{M}^{(v)}_2$ respectively. If
$$\left ((A_j)_{1\leq j\leq n},(G^{(v)}_i)_{\overset{1\leq i\leq
r_v}{ 1\leq v\leq f}}, (H^{(v)}_l)_{\overset{1 \leq l\leq
s_v}{1\leq v\leq f}}, \lbrace E^{(v)}_{jl} \rbrace _{j,l,v},
\lbrace F^{(v)}_{jl} \rbrace _{j,l,v},(Z_u)_{1\leq u\leq
 p}\right )$$ is an arbitrary microstate in the set of matricial microstates
\begin{eqnarray}&&\Gamma_R\bigg((x_j)_{1\leq j\leq n},(p^{(v)}_i)_{\overset{1\leq
i\leq r_v}{ 1\leq v\leq f}}, (q^{(v)}_l)_{\overset{1\leq l \leq
s_v}{1\leq v\leq f}}, \lbrace e^{(v)}_{jl} \rbrace
_{j,l,v},\lbrace f^{(v)}_{jl} \rbrace
_{j,l,v},\\\nonumber&&\hspace{1 cm}(z_u)_{1\leq u\leq
p};m,k,\epsilon\bigg)
\end{eqnarray}
and if $m$ is large and $\epsilon$ is small enough, then
$$\left|\left|A_j-\phi _j\left((G^{(v)}_i)_{\overset{1\leq i\leq r_v}{ 1\leq v\leq f}},
(H^{(v)}_l)_{\overset{1 \leq l\leq s_v}{ 1\leq v\leq
f}},(Z_u)_{1\leq u\leq
 p}\right)\right|\right|_2<\omega,\,\, j=1,\ldots,n.$$
Let $\delta >0$ and write $k=k_0^2 t+w$ for some integers $w$, $t$
with $0\leq w\leq k_0^2 -1$. If $m,\epsilon$ are suitably chosen,
then there exist $\mc{M}^{(v)}_1\cong
\tilde{\mc{M}}^{(v)}_1\subset \mc{M}_k(\mb{C})$,
$\mc{M}^{(v)}_2\cong \tilde{\mc{M}}^{(v)}_2\subset
\mc{M}_k(\mb{C})$ (not necessarily unital inclusions) and matrix
units $\lbrace \tilde{E}^{(v)}_{jl} \rbrace _{j,l,v}\subset
\tilde{\mc{M}}^{(v)}_1$, $\lbrace \tilde{F}^{(v)}_{jl} \rbrace
_{j,l,v}\subset \tilde{\mc{M}}^{(v)}_2$, such that
$$\left|\left|\tilde{E}^{(v)}_{jl}-E^{(v)}_{jl}\right|\right|_2<\delta,
\,\,\left|\left|\tilde{F}^{(v)}_{jl}-F^{(v)}_{jl}\right|\right|_2<
\delta \,\,\forall\,1\leq j ,l\leq k_0,$$ and
$\tilde{\mc{M}}^{(v)}_1\subset\left(\tilde{\mc{M}}^{(v)}_2\right)'\cap
\mc{M}_k(\mb{C})$. The relative commutants of
$\tilde{\mc{M}}^{(v)}_1$ and $\tilde{\mc{M}}^{(v)}_2$ in
$\mc{M}_k(\mb{C})$ satisfy
$\left(\tilde{\mc{M}}^{(v)}_1\right)'\cap \mc{M}_k(\mb{C})\cong
(\mc{M}_{k_0}(\mb{C})\otimes 1\otimes \mc{M}_t(\mb{C}))\oplus
\mc{M}_w(\mb{ C})$ and $\left(\tilde{\mc{M}}^{(v)}_2\right)'\cap
\mc{M}_k(\mb{C})\cong (1\otimes \mc{M}_{k_0}(\mb{C})\otimes
\mc{M}_t(\mb{C}))\oplus \mc{M}_w(\mb{ C})$. Let $\eta^{(v)}
(x,\lbrace e^{(v)}_{jl}\rbrace _{j,l}):=
\frac{1}{k_0}\sum_{j,l=1}^{k_0}e^{(v)}_{jl}xe^{(v)}_{lj} \in
\mb{C}<X_1,\ldots,X_{k_0^2 +1}>$ be the polynomial in $k_0^2+1$
indeterminates that gives the conditional expectation
$E_{(\mc{M}^{(v)}_1)' \cap \mc{M}}: \mc{M}\rightarrow
(\mc{M}^{(v)}_1)' \cap \mc{M}$, $E_{(\mc{M}^{(v)}_1)' \cap
\mc{M}}(x)= \eta^{(v)} (x,\lbrace e^{(v)}_{jl} \rbrace _{j,l})$.
Then $G^{(v,1)}_1:=\eta^{(v)} (G^{(v)}_1,$ $
\lbrace\tilde{E}^{(v)}_{jl}\rbrace _{j,l}) \in
\left(\tilde{\mc{M}}^{(v)}_1\right)'\cap \mc{M}_k(\mb{C})$ and
since $p^{(v)}_1=E_{(\mc{M}^{(v)}_1)' \cap
\mc{M}}(p^{(v)}_1)=\eta^{(v)} (p^{(v)}_1,$ $\lbrace
e^{(v)}_{jl}\rbrace _{j,l})$ it follows that
$$\left|\tau_k \left((G^{(v,1)}_1)^l\right)-\tau \left((p^{(v)}_1)^l\right)\right|<
\delta_1,\,\,\forall 1\leq l \leq m_1$$ for any given
$\delta_1,m_1$, provided that $\epsilon ,\delta$ are small and $m$
is large enough. For suitable $m_1$, $\delta_1$ there exists
 a projection $P^{(v,1)}_1 \in \left(\tilde{\mc{M}}^{(v)}_1\right)'\cap \mc{M}_k(\mb{C})$
 of rank
$\left[\frac{k_0 t+w}{2}\right]$ such that
$||P^{(v,1)}_1-G^{(v,1)}_1||_2<\delta_2$. Then
$||G^{(v)}_1-P^{(v,1)}_1||_2\leq
||G^{(v)}_1-G^{(v,1)}_1||_2+||G^{(v,1)}_1-P^{(v,1)}_1||_2
<2\delta_2$ since $||G^{(v)}_1-G^{(v,1)}_1||_2<\delta_2$ for
convenient $m$, $\epsilon$, $\delta$. With this procedure we can
find projections $P^{(v,1)}_1,\ldots,P^{(v,1)}_{r_v}\in
\left(\tilde{\mc{M}}^{(v)}_1\right)' \cap \mc{M}_k(\mb{C})$ and
$Q^{(v,1)}_1,\ldots ,Q^{(v,1)}_{s_v}\in
\left(\tilde{\mc{M}}^{(v)}_2\right)'\cap \mc{M}_k(\mb{C})$, all of
rank $\left[\frac{k_0 t+w}{2}\right]$, such that
$||G^{(v)}_i-P^{(v,1)}_i||_2<2 \delta_2$ and
$||H^{(v)}_j-Q^{(v,1)}_j||_2<2\delta_2$ for all indices $i,j,v$.
Moreover,
$$\left|\left|A_j-\phi_j\left((P^{(v,1)}_i)_{\overset{1\leq i\leq r_v}{ 1\leq v\leq f}},
(Q^{(v,1)}_l)_{\overset{1\leq l\leq s_v}{1\leq v\leq
f}},(Z_u)_{1\leq u\leq
 p}\right)\right|\right|_2<\omega\,\forall\,1\leq j\leq n$$
if we choose a sufficiently small $\delta_2>0$. Let
$\mc{G}^{(v)}_1(k)\subset \left(\tilde{\mc{M}}^{(v)}_1\right)'\cap
\mc{M}_k(\mb{C})$ and
$\mc{G}^{(v)}_2(k)\subset\left(\tilde{\mc{M}}^{(v)}_2\right)'\cap
\mc{M}_k(\mb{C})$ be $2$ fixed copies of the Grassmann manifold
$\mc{G}\left( k_0 t+w,\left[\frac{k_0t+w}{2}\right]\right)$
(projections in $\mc{M}_{k_0t+w}(\mb{C})$, of rank
$\left[\frac{k_0t+w}{2}\right]$). There exists a unitary
$U^{(v)}\in \mc{U}(k)$ such that $U^{(v)}P^{(v,1)}_1U^{(v)*},$
$\ldots ,U^{(v)}P^{(v,1)}_{r_v}U^{(v)*}\in \mc{G}^{(v)}_1(k)$ and
$U^{(v)}Q^{(v,1)}_1U^{(v)*},\ldots,U^{(v)}Q^{(v,1)}_{s_v}U^{(v)*}\in
\mc{G}^{(v)}_2(k)$. The previous inequality becomes
\begin{eqnarray}&&\bigg|\bigg| A_j-\phi_j\bigg( (U^{(v)}P^{(v,1)}_iU^{(v)*})_{\overset{1\leq
i\leq r_v}{1\leq v\leq f}},
(U^{(v)}Q^{(v,1)}_lU^{(v)*})_{\overset{1\leq l\leq s_v}{ 1\leq
v\leq f}},(Z_u)_{1\leq u\leq p},\\\nonumber&&\hspace{1 cm}
(\mbox{Re}(U^{(v)}),\mbox{Im}(U^{(v)}))_{1\leq v\leq
f}\bigg)\bigg|\bigg|_2<\omega\,\forall\, 1\leq j\leq n.
\end{eqnarray}
The euclidean norm on $\mc{M}_k^{sa}$ induces a $\mc{U}(k_0
t+w)$-invariant metric on the manifold $\mc{G}\left( k_0
t+w,\left[\frac{k_0t+w}{2}\right]\right)$ and if $\lbrace P_a
\rbrace _{a\in A(k)}$ is a minimal $\theta$-net in the manifold
with respect to this metric, then (\cite {12})
 $|A(k)|\leq \left(\frac{Ch_k}{\theta}\right)^{g_k}$ where $C$ is a
 universal constant, $g_k=2\left[\frac{k_0t+w}{2}\right]\cdot
 \left(k_0 t+w-\left[\frac{k_0t+w}{2}\right]\right)$ is the dimension of
 $\mc{G}\left( k_0 t+w,\left[\frac{k_0t+w}{2}\right]\right)$ and
 $h_k\leq \sqrt {2k}$ is the diameter of the Grassmann manifold
 $\mc{G}\left( k_0 t+ w,\left[\frac{k_0t+w}{2}\right]\right)$ in $\mc{M}_k^{sa}.$
 There exist $\alpha :=(a^{(v)}_1,\ldots,a^{(v)}_{r_v})_{1\leq v\leq f}$
 and $\beta :=(b^{(v)}_1,\ldots,b^{(v)}_{s_v})_{1\leq v\leq f}$
 with entries from $A(k)$
 such that $$\left|\left|P^{(v)}_{a^{(v)}_i}-U^{(v)}P^{(v,1)}_iU^{(v)*}\right|\right|_e\leq
\theta\,\,\mbox{and}\,\,\left|\left|P^{(v)}_{b^{(v)}_l}
-U^{(v)}Q^{(v,1)}_lU^{(v)*}\right|\right|_e\leq\theta$$ for all
$1\leq i\leq r_v$, $1\leq l\leq s_v$, $1\leq v\leq f$. The
polynomials $(\phi_j)_{1\leq j\leq n}$ are in particular Lipschitz
functions hence there exists a constant $D=D\left((\phi_j)_{1\leq
j\leq n}, R\right)>0$ (note that $|\alpha |= r_1+\ldots +r_f$ and
$|\beta |=s_1+\ldots +s_f$) such that
\begin{eqnarray}&&
\left|\left|\phi_j(V_1,\ldots,V_{|\alpha|+|\beta|+p+2f})-\phi_j(W_1,
\ldots,W_{|\alpha|+|\beta|+p+2f})\right|\right|_e\\\nonumber&&\hspace{1
cm} \leq D\left|\left|(V_1,\ldots,V_{|\alpha|+|\beta|+p+2f})-
(W_1,\ldots,W_{|\alpha|+|\beta|+p+2f})\right|\right|_e
\end{eqnarray}
for all $1\leq j\leq n$ and all
$V_1,\ldots,V_{|\alpha|+|\beta|+p+2f},W_1,
\ldots,W_{|\alpha|+|\beta|+p+2f} \in\lbrace V\in
M_k\,\,|\,\,||V||\leq R\rbrace$. We have then
\begin{eqnarray}&&\bigg|\bigg|A_j-\phi_j\bigg((P_a)_{a\in\alpha},(P_b)_{b\in\beta},
(Z_u)_{1\leq u\leq p},
\big(\mbox{Re}(U^{(v)}),\mbox{Im}(U^{(v)})\big)_{1\leq v\leq
f}\bigg)\bigg|\bigg|_e\\\nonumber&&\hspace{1 cm}
<\omega\sqrt{k}+D\bigg|\bigg|\bigg((U^{(v)}P^{(v,1)}_iU^{(v)*})_{\overset{1\leq
i\leq r_v}{1\leq v\leq f}},
(U^{(v)}Q^{(v,1)}_lU^{(v)*})_{\overset{1\leq l\leq s_v}{1\leq
v\leq f}},\\\nonumber&&\hspace{1 cm}(Z_u)_{1\leq u\leq
p},\big(\mbox{Re}(U^{(v)}),\mbox{Im}(U^{(v)})\big)_{1\leq v\leq
f}\bigg)-
\bigg((P_a)_{a\in\alpha},(P_b)_{b\in\beta},\\\nonumber&&\hspace{1
cm} (Z_u)_{1\leq u\leq p},\big(\mbox{Re}(U^{(v)}),
\mbox{Im}(U^{(v)})\big)_{1\leq v\leq
f}\bigg)\bigg|\bigg|_e<\omega\sqrt{k}\\\nonumber&& \hspace{1
cm}+D\theta\sqrt{|\alpha|+|\beta|}=2\omega\sqrt{k},
\end{eqnarray}
if we choose $\theta
:=\frac{\omega}{D}\sqrt{\frac{k}{|\alpha|+|\beta|}}$. Define
$F_{\alpha ,\beta}:(\mc{M}_k^{sa})^{p+2f}\rightarrow
(\mc{M}_k^{sa})^n$ by
\begin{eqnarray}&&F_{\alpha
,\beta}\left((W_u)_{1\leq u\leq p},(W^{(v)}_1,W^{(v)}_2)_{1\leq v
\leq f}\right)\\\nonumber&&\hspace{1 cm} =\left(\phi_j
((P_a)_{a\in\alpha},(P_b)_{b\in\beta},(W_u)_{1\leq u\leq p},
(W^{(v)}_1,W^{(v)}_2)_{1\leq v\leq f})\right)_{1\leq j\leq n},
\end{eqnarray}
and note that $\mbox{dist}_e((A_j)_{1\leq j\leq
n},\mbox{Range}(F_{\alpha ,\beta})) <2\omega\sqrt{nk}$. Note also
that all the components of $F_{\alpha ,\beta}$ are polynomial
functions of degrees $\leq 3d+2$. Use now Lemma \ref{lema1} as in
the proof of Proposition \ref{prop1} to obtain the estimates:
\begin{eqnarray}&&
\mbox{vol}_{nk^2}\bigg(\Gamma_R\bigg((x_j)_{1\leq j\leq
n}:(p^{(v)}_i)_{\overset{1\leq i\leq r_v}{1\leq v\leq f}},
(q^{(v)}_l)_{\overset{1\leq l \leq s_v}{1\leq v\leq f}},\lbrace
e^{(v)}_{jl} \rbrace _{j,l,v}, \lbrace f^{(v)}_{jl} \rbrace
_{j,l,v},\\\nonumber&&\hspace{.5 cm} (z_u)_{1\leq u\leq
p};m,k,\epsilon\bigg)\bigg) \leq
\left(\left(\frac{Ch_k}{\theta}\right)^{g_k}\right)^{|\alpha|+|\beta|}\cdot
\binom{nk^2} {(p+2f)k^2}\cdot C(d)\\\nonumber&&\hspace{.5 cm}
\cdot \mbox{vol}_{(p+2f)k^2}\left(B(0,(a+2\omega)\sqrt{nk})\right)
\cdot \mbox{vol}_{nk^2-(p+2f)k^2}\left(B(0,2\omega
\sqrt{nk})\right)\\\nonumber&&\hspace{.5 cm}=
\left(\frac{CDh_k}{\omega}\sqrt{\frac{|\alpha|+|\beta|}{k}}
\right)^{(|\alpha|+|\beta|)g_k} \cdot \binom{nk^2}{(p+2f)k^2}
\cdot C(d)\\\nonumber&&\hspace{.5 cm} \cdot \frac{(\pi
nk)^{\frac{(p+2f)k^2}{2}}(2\omega +a)^{(p+2f)k^2}}{\Gamma
\left(1+\frac{(p+2f)k^2}{2}\right)}\cdot \frac{(\pi
nk)^{\frac{nk^2-(p+2f)k^2}{2}}(2\omega)^{nk^2-(p+2f)k^2}}{\Gamma
\left(1+\frac{nk^2-(p+2f)k^2}{2}\right)}.
\end{eqnarray}
The above estimate, the inequality (\ref{g}) on pag. \pageref{g},
and the inequalities
\begin{eqnarray}&&h_k\leq
\sqrt{2k},0<\omega\leq a, g_k=2\left[\frac{k_0 t+w}{2}\right]
\left(k_0 t+w-\left[\frac{k_0 t+w}{2}\right]\right)\leq
2\\\nonumber&&\hspace{.5 cm}\cdot\frac{k_0 t+w}{2}\cdot\left(k_0
t+ w-\frac{k_0 t+w}{2}\right)=\frac{(k_0 t+w)^2} {2}=\frac{(k+k_0
w-w)^2}{2k_0^2},\end{eqnarray} together with $C(d)\leq
(3d+2)^{(p+2f)k^2}$ imply
\begin{eqnarray}&&\mbox{vol}_{nk^2}\bigg(\Gamma_R\bigg((x_j)_{1\leq j\leq
n}:(p^{(v)}_i)_{\overset{1\leq i\leq r_v}{1\leq v\leq f}},
(q^{(v)}_l)_{\overset{1\leq l \leq s_v}{1\leq v\leq f}},\lbrace
e^{(v)}_{jl} \rbrace _{j,l,v}, \lbrace f^{(v)}_{jl} \rbrace
_{j,l,v},\\\nonumber&&\hspace{.5 cm}(z_u)_{1\leq u\leq
p};m,k,\epsilon\bigg)\bigg)\leq
\left(\frac{CD\sqrt{2(|\alpha|+|\beta|)}}{\omega}\right)^{\frac
{(k+k_0 w-w)^2} {2k_0^2}(|\alpha|+|\beta|)}\\\nonumber&&\hspace{.5
cm} \cdot \frac{2^{\frac{nk^2}{2}}(\pi
nk)^{\frac{nk^2}{2}}(3a)^{(p+2f)k^2}(2\omega)^{(n-p-2f)k^2}}
{\Gamma \left(1+\frac{nk^2}{2}\right)}\\\nonumber&&\hspace{.5 cm}
\cdot\binom{nk^2}{(p+2f)k^2}\cdot (3d+2)^{(p+2f)k^2}
\end{eqnarray}
therefore \begin{eqnarray}&&\frac{1}{k^2}\chi_R \bigg((x_j)_{1\leq
j\leq n}:(p^{(v)}_i)_{\overset{1\leq i\leq r_v}{1\leq v\leq f}},
(q^{(v)}_l)_{\overset{1\leq l \leq s_v}{1\leq v\leq f}}, \lbrace
e^{(v)}_{jl} \rbrace _{j,l,v}, \lbrace f^{(v)}_{jl} \rbrace
_{j,l,v},\\\nonumber&&\hspace{.5 cm} (z_u)_{1\leq u\leq p};
m,k,\epsilon\bigg) +\frac{n}{2}\log k\leq C(n,p,a,d,f)+n\log
k\\\nonumber&&\hspace{.5 cm} +\frac{|\alpha|+|\beta|}
{2k_0^2}\left(1+\frac{k_0w-w}{k}\right)^2\log\frac{CD\sqrt{2(|\alpha|+|\beta|)}}
{\omega}\\\nonumber&&\hspace{.5 cm}+(n-p-2f)\log\omega
-\frac{1}{k^2}\log\Gamma\left(1+\frac{nk^2}{2}\right)+
\frac{1}{k^2}\log \binom{nk^2}{(p+2f)k^2}.
\end{eqnarray}
Use $\frac{1}{k^2}\log \binom{nk^2}{(p+2f)k^2}=n\log n-(p+2f)\log
(p+2f)-(n-p-2f)\log (n-p-2f)+o(1)$ and Stirling's formula
$\frac{1}{k^2}\log\Gamma\left(1+\frac{nk^2}{2}\right)=\frac{n}{2}
\log\frac{nk^2}{2e}+o(1)$ to conclude
\begin{eqnarray}&&\label{estimate2}
\chi_R\bigg((x_j)_{1\leq j\leq n}:(p^{(v)}_i)_{\overset{1\leq
i\leq r_v}{1\leq v\leq f}}, (q^{(v)}_l)_{\overset{1\leq l \leq
s_v}{1\leq v\leq f}}, \lbrace e^{(v)}_{jl} \rbrace _{j,l,v},
\lbrace f^{(v)}_{jl} \rbrace _{j,l,v},\\\nonumber&&\hspace{.5 cm}
(z_u)_{1\leq u\leq p};m,\epsilon\bigg) \leq
\frac{|\alpha|+|\beta|}{2k_0^2}\log(CD\sqrt{2(|\alpha|+|\beta|)})\\
\nonumber&&\hspace{.5 cm} +C(n,p,a,d,f)+\left(n-p-2f-\frac{
|\alpha|+|\beta|}{2k_0^2}\right)\log\omega.
\end{eqnarray}
The last inequality shows that the free entropy of $\{x_1,\ldots ,
x_n\}$ does not exceed $C(n,p,a,d,f)+(n-p-2f)\log\omega$ since
$k_0$ is an arbitrary integer, $R$ is an arbitrary positive number
and $x_1,\ldots , x_n$ generate $M$.
\end{proof}
\subsection{Hyperfinite dimension of free group factors}\label{hyp}
\begin{thm}\label{prop4}
If $n\geq p+2f+1$, then the free group factor $\mc{L}(\mb{F}_n)$
can not be asymptotically decomposed as
$$\lim_{\omega\rightarrow 0}\,^{||\cdot ||_2}\sum_{\overset{1\leq j_1,\ldots ,j_{t+1}\leq
f}{1\leq t\leq d}} \mc{N}_{j_1}^\omega \mc{Z}^\omega
\mc{N}_{j_2}^\omega \mc{Z}^\omega \ldots \mc{N}_{j_t}^\omega
\mc{Z}^\omega \mc{N}_{j_{t+1}}^\omega$$  where
$\{\mc{Z}^\omega\subset \mc{L}(\mb{F}_n)\}_\omega$ contain $p$
self-adjoint elements, $\{\mc{N}_1^\omega,\ldots
,\mc{N}_f^\omega\}_\omega$ are nonprime subfactors of
$\mc{L}(\mb{F}_n)$, and $d\geq 1$ is an integer.
\end{thm}
\begin{proof} Suppose first that $\infty >n\geq p+2f+1$ and consider a semicircular
system $\{x_1,\ldots ,x_n\}$ that generates $\mc{L}(\mb{F}_n)$ as
a von Neumann algebra. If the assertion were true then one could
find for every $\omega >0$ noncommutative polynomials and
projections as in Proposition \ref{prop3}, satisfying the
inequalities (\ref{cond1}). But then the estimate of the free
entropy (\ref{estimate1}) would imply that
$\chi(x_1,\ldots,x_n)=-\infty$ if one makes $\omega$ tend to $0$,
contradiction.

If $n=\infty$ then $\mc{L}(\mb{F}_\infty)$ is generated by an
infinite semicircular system $\{x_t\}_{t\geq 1}$. If we fix an
integer $k\geq p+2f+1$, then we can approximate $x_1,\ldots ,x_k$
by polynomials $(\phi_j)_{1\leq j\leq k}$ as in (\ref{cond1}) and
so one has the estimate of the modified free entropy
(\ref{estimate2}) with $k$ instead of $n$. Taking $m$,
$\frac{1}{\epsilon}$, $R$, $k_0\rightarrow\infty$ and
$\omega\rightarrow 0$ in this estimate, one obtains
\begin{eqnarray}&&\chi\bigg((x_j)_{1\leq j\leq k}:(p^{(v)}_i)_{\overset{1\leq
i\leq r_v}{1\leq v\leq f}}, (q^{(v)}_l)_{\overset{1\leq l \leq
s_v}{1\leq v\leq f}}, \lbrace e^{(v)}_{jl} \rbrace _{j,l,v},
\lbrace f^{(v)}_{jl} \rbrace _{j,l,v},\\\nonumber&&\hspace{1 cm}
(z_u)_{1\leq u\leq p}\bigg) <\chi (x_1,\ldots ,x_k)\end{eqnarray}
where $(p^{(v)}_i)_{\overset{1\leq i\leq r_v}{1\leq v\leq f}}$,
$(q^{(v)}_l)_{\overset{1\leq l \leq s_v}{1\leq v\leq f}}$,
$\lbrace e^{(v)}_{jl} \rbrace _{j,l,v}$, $\lbrace f^{(v)}_{jl}
\rbrace _{j,l,v}$, $(z_u)_{1\leq u\leq p}$ are as in Proposition
\ref{prop3}. If $\mc{A}_t$ denotes the von Neumann algebra
$\{x_1,\ldots ,x_t\}''$ and $E_t$ the conditional expectation onto
it, then \begin{eqnarray}&&\bigg((x_j)_{1\leq j\leq
k},(E_t(p^{(v)}_i))_{\overset{1\leq i\leq r_v} {1\leq v\leq
f}},(E_t(q^{(v)}_l))_{\overset{1\leq l \leq s_v}{ 1\leq v\leq
f}},\\\nonumber&&\hspace{1 cm} \lbrace E_t(e^{(v)}_{jl}) \rbrace
_{j,l,v}, \lbrace E_t(f^{(v)}_{jl}) \rbrace
_{j,l,v},(E_t(z_u))_{1\leq u\leq p}\bigg)_{t\geq 1}\end{eqnarray}
converges in distribution as $t\rightarrow\infty$ to
$$\left((x_j)_{1\leq j\leq k},(p^{(v)}_i)_{\overset{1\leq i\leq r_v}{
1\leq v\leq f}}, (q^{(v)}_l)_{\overset{1\leq l \leq s_v}{ 1\leq
v\leq f}}, \lbrace e^{(v)}_{jl} \rbrace _{j,l,v}, \lbrace
f^{(v)}_{jl} \rbrace _{j,l,v}, (z_u)_{1\leq u\leq p}\right)$$
therefore \begin{eqnarray}&&\chi\bigg((x_j)_{1\leq j\leq
k}:(E_t(p^{(v)}_i))_{\overset{1\leq i\leq r_v} {1\leq v\leq f}},
(E_t(q^{(v)}_l))_{\overset{1\leq l \leq s_v}{ 1\leq v\leq f}},
\lbrace E_t(e^{(v)}_{jl}) \rbrace _{j,l,v},\\\nonumber&&\hspace{1
cm} \lbrace E_t(f^{(v)}_{jl}) \rbrace _{j,l,v},(E_t(z_u))_{1\leq
u\leq p}\bigg) <\chi (x_1,\ldots ,x_k)\end{eqnarray} for some
large integer $t>k$. But this leads to a contradiction:
\begin{eqnarray}&&\chi (x_1,\ldots ,x_t)= \chi\bigg((x_j)_{1\leq j\leq
t}:(E_t(p^{(v)}_i))_{\overset{1\leq i\leq r_v}{ 1\leq v\leq f}},
(E_t(q^{(v)}_l))_{\overset{1\leq l \leq s_v}{ 1\leq v\leq
f}},\\\nonumber&&\hspace{.5 cm} \lbrace E_t(e^{(v)}_{jl}) \rbrace
_{j,l,v}, \lbrace E_t(f^{(v)}_{jl}) \rbrace _{j,l,v},
(E_t(z_u))_{1\leq u\leq p}\bigg)\leq \chi\bigg((x_j)_{1\leq j\leq
k}:\\\nonumber&&\hspace{.5 cm} (E_t(p^{(v)}_i))_{\overset{1\leq
i\leq r_v}{ 1\leq v\leq f}}, (E_t(q^{(v)}_l))_{\overset{1\leq l
\leq s_v}{ 1\leq v\leq f}}, \lbrace E_t(e^{(v)}_{jl}) \rbrace
_{j,l,v}, \lbrace E_t(f^{(v)}_{jl}) \rbrace
_{j,l,v},\\\nonumber&&\hspace{.5 cm}(E_t(z_u))_{1\leq u\leq
p}\bigg)+\chi (x_{k+1},\ldots ,x_t)<\chi (x_1,\ldots
,x_k)\\\nonumber&&\hspace{.5 cm}+\chi (x_{k+1},\ldots ,x_t) =\chi
(x_1,\ldots ,x_t).
\end{eqnarray}
\end{proof}
\begin{cor}\label{cor6}
If $\mc{P}\subset \mc{L}(\mb{F}_n)$ is a subfactor of finite index
and if $n\geq p+2f+2$, then $\mc{P}$ can not be asymptotically
decomposed as
$$\lim_{\omega\rightarrow 0}\,^{||\cdot ||_2}\sum_{\overset{1\leq j_1,
\ldots ,j_{t+1}\leq f}{ 1\leq t\leq d}} \mc{N}_{j_1}^\omega
\mc{Z}^\omega \mc{N}_{j_2}^\omega \mc{Z}^\omega \ldots
\mc{N}_{j_t}^\omega \mc{Z}^\omega \mc{N}_{j_{t+1}}^\omega ,$$
where $\{\mc{Z}^\omega\}_\omega$ contain $p$ self-adjoint elements
of $\mc{P}$, $\{\mc{N}_1^\omega,\ldots , \mc{N}_f^\omega\}_\omega$
are nonprime subfactors of $\mc{P}$, and $d\geq 1$ is an integer.
\end{cor}
\begin{proof} Since $\mc{P}\subset \mc{L}(\mb{F}_n)$ is a subfactor of finite index,
$\mc{L}(\mb{F}_n)$ can be obtained from $\mc{P}$ with the basic
construction (\cite{7}, \cite{8}): there exists a subfactor
$\mc{Q}\subset \mc{P}$ such that
$\mc{L}(\mb{F}_n)=<\mc{P},e_\mc{Q}>$, where $e_\mc{Q}$ is the
Jones projection associated to the inclusion $\mc{Q}\subset
\mc{P}$. But $<\mc{P},e_\mc{Q}>=\mc{P} e_\mc{Q}\mc{P}$ (\cite{8}),
hence $\mc{L}(\mb{F}_n)$ can be decomposed as
$\mc{P}e_\mc{Q}\mc{P}$. Apply now Theorem \ref{prop4}.
\end{proof}
\begin{cor}\label{chn}
If $n\geq p+2f+1$, then the free group factor $\mc{L}(\mb{F}_n)$
can not be decomposed as
$$\overline{\s}^w\sum_{\overset{1\leq j_1,\ldots ,j_{t+1}\leq
f}{ 1\leq t\leq d}} \mc{N}_{j_1}\mc{Z}\mc{N}_{j_2}\mc{Z}\ldots
\mc{N}_{j_t}\mc{Z}\mc{N}_{j_{t+1}} ,$$  where $\mc{Z}\subset
\mc{L}(\mb{F}_n)$ contains $p$ self-adjoint elements,
$\mc{N}_1,\ldots ,\mc{N}_f$ are nonprime subfactors of
$\mc{L}(\mb{F}_n)$, and $d\geq 1$ is an integer. Moreover, if
$\mc{P}\subset \mc{L}(\mb{F}_n)$ is a subfactor of finite index
and if $n\geq p+2f+2$, then $\mc{P}$ also can not be decomposed as
$$\overline{\s}^w\sum_{\overset{1\leq j_1,\ldots ,j_{t+1}\leq
f}{ 1\leq t\leq d}} \mc{N}_{j_1} \mc{Z}\mc{N}_{j_2}\mc{Z} \ldots
\mc{N}_{j_t} \mc{Z}\mc{ N}_{j_{t+1}} ,$$ for any subset $\mc{Z}$
containing $p$ self-adjoint elements of $\mc{P}$, any
$\mc{N}_1,\ldots ,$ $\mc{N}_f$ nonprime subfactors of $\mc{P}$,
and any integer $d\geq 1$.
\end{cor}
\begin{proof} Follows from Theorem \ref{prop4} and
Corollary \ref{cor6}, for $\mc{Z}^\omega =\mc{Z}$,
$\mc{N}_1^\omega =\mc{N}_1,\ldots, \mc{N}_f^\omega =\mc{N}_f$.
\end{proof}
Corollary \ref{cor3} settles a conjecture of L. Ge and S. Popa
(\cite{5}) in the case $n=\infty$. Recall that for a type
$\mbox{\!I\!I}_1$-factor $\mc{M}$ one defines
$\ell_h(\mc{M})=\min\{f\in\mb{N}\,\,|\,\,\exists\,\,
\mbox{hyperfinite}\,\,\mc{R}_1,\ldots ,\mc{R}_f \subset
M\,\,\mbox{s.t.}\,\, \overline{\mbox{sp}}^w
\mc{R}_1\mc{R}_2\ldots$ $\mc{R}_f=\mc{M}\}$. Note that the
definition of hyperfinite dimension is given in terms of
hyperfinite subalgebras. If one defined the hyperfinite dimension
in terms of hyperfinite subfactors instead of hyperfinite
subalgebras, then the proof of Corollary \ref{cor3} would have
followed immediately from Corollary \ref{chn}. But with Definition
\ref{hd}, we need the asymptotic indecomposability result from
Theorem \ref{prop4}.
\begin{cor}\label{cor3}
$\ell_h(\mc{L}(\mb{F}_n))\geq [\frac{n-2}{2}]+1\,\,\forall 4\leq
n\leq\infty$.
\end{cor}
\begin{proof} If $\ell_h(\mc{L}(\mb{F}_n))\leq [\frac{n-2}{2}]$, then $\mc{L}(\mb{
F}_n)=\overline{\mbox{sp}}^w\mc{R}_1\mc{R}_2\ldots \mc{R}_f$ for
some hyperfinite subalgebras $\mc{R}_1,\ldots ,\mc{R}_f$ and some
integer $f$ with $n\geq 2f+2$. Let $m\geq 1$, $y_1,\ldots ,y_m\in
\mc{L}(\mb{F}_n)$ and $\omega >0$ be fixed. Then there exist
finite dimensional subalgebras $\mc{B}_v^\omega
=\mc{B}_v(y_1,\ldots ,y_m;\omega )\subset \mc{R}_v$, $1\leq v\leq
f$, such that $$\mbox{dist}_{||\cdot
||_2}\left(y_j,\mc{B}_1^\omega \mc{B}_2^\omega\ldots
\mc{B}_f^\omega\right)<\omega\,\,\forall 1\leq j\leq m.$$ Each
finite dimensional subalgebra $\mc{B}_v^\omega$ is contained in a
copy of the hyperfinite $\mbox{\!I\!I}_1$-factor, say
$\mc{B}_v^\omega\subset \mc{R}_v^\omega =\mc{R}_v^\omega
(y_1,\ldots ,y_m;\omega )\subset \mc{L}(\mb{ F}_n)$. Consequently,
$$\mbox{dist}_{||\cdot ||_2}\left(y_j,\mc{R}_1^\omega \mc{R}_2^\omega \ldots
\mc{R}_f^\omega\right)<\omega\,\, \forall 1\leq j\leq m,$$ hence
$\mc{L}(\mb{F}_n)$ admits an asymptotic decomposition of the form
$$\lim_{\omega\rightarrow 0}\,^{||\cdot ||_2}\mc{R}_1^\omega
\mc{R}_2^\omega\ldots \mc{R}_f^\omega,$$ in contradiction with
Theorem \ref{prop4} as $\mc{R}_1^\omega ,\ldots ,\mc{R}_f^\omega$
are nonprime and $n\geq 2f+2$.
\end{proof}
\begin{cor}\label{cor7}
If $\mc{P}\subset \mc{L}(\mb{F}_n)$ is a subfactor of finite index
and $5\leq n\leq\infty$, then $\ell_h(\mc{P})\geq
[\frac{n-3}{2}]+1$.
\end{cor}
\begin{proof} Follows from Corollary \ref{cor6}.
\end{proof}
\section{Indecomposability over abelian subalgebras}\label{chab}
\setcounter{equation}{0} Another estimate of free entropy is used
to prove that the free group factor $\mc{L}(\mb{F}_n)$ does not
admit an asymptotic decomposition of the form
$$\lim_{\omega\rightarrow 0}\,^{||\cdot ||_2}\sum_{\overset{1\leq j_1,\ldots
,j_{t+1}\leq f}{ 1\leq t\leq d}} \mc{A}_{j_1}^\omega \mc{Z}^\omega
\mc{A}_{j_2}^\omega \mc{Z}^\omega\ldots \mc{A}_{j_t}^\omega
\mc{Z}^\omega \mc{A}_{j_{t+1}}^\omega ,$$ where $\{\mc{A}_1^\omega
,\ldots ,\mc{A}_f^\omega\}$ are abelian subalgebras of
$\mc{L}(\mb{F}_n)$, $\{\mc{Z}^\omega\subset
\mc{L}(\mb{F}_n)\}_\omega$ are subsets with $p$ self-adjoint
elements, $d\geq 1$ is an arbitrary integer, and $n \geq p+2f +1$.
Similarly, for free group subfactors one has the following: if
$n\geq p+2f+2$ and $\mc{P}\subset \mc{L}(\mb{F}_n)$ is a subfactor
of finite index, then $\mc{P}$ does not admit such an asymptotic
decomposition either. In particular, the abelian dimension of
$\mc{L}(\mb{F}_n)$ is $\geq [\frac{n-2}{2}]+1$ and the abelian
dimension of $\mc{P}$ is $\geq [\frac{n-3}{2}]+1$. For $n=\infty$
this proves the second part of L. Ge's and S. Popa's (\cite{5})
conjecture: the abelian dimension of free group factors is
infinite. The definitions of abelian dimension and asymptotic
decomposition over abelian subalgebras are given next.
\begin{defn}
(\cite{5}) If $\mc{M}$ is a $\mbox{\!I\!I}_1$-factor, then the
abelian dimension of $\mc{M}$, denoted $\ell_a(\mc{M})$, is
defined as the smallest positive integer $f\in\mb{N}$ with the
property that there exist abelian subalgebras $\mc{A}_1,\ldots
,\mc{A}_f\subset \mc{M}$ such that
$\overline{\s}^w\mc{A}_1\mc{A}_2\ldots \mc{A}_f=\mc{M}$. If there
is no such positive integer $f$, then by definition,
$\ell_a(\mc{M})=+\infty$.
\end{defn}
\begin{defn}\label{ade}
A type $\mbox{\!I\!I}_1$-factor $\mc{M}$ admits an asymptotic
decomposition over abelian subalgebras, denoted
$$\lim_{\omega\rightarrow 0}\,^{||\cdot ||_2}\sum_{\overset{1\leq j_1,\ldots
,j_{t+1}\leq f}{1\leq t\leq d}} \mc{A}_{j_1}^\omega \mc{Z}^\omega
\mc{A}_{j_2}^\omega\mc{Z}^\omega\ldots \mc{A}_{j_t}^\omega
\mc{Z}^\omega \mc{A}_{j_{t+1}}^\omega,$$ provided that $\forall
n\geq 1$ $\forall x_1,\ldots ,x_n\in \mc{M}$ $\forall\omega >0$
$\exists \mc{A}_1^\omega = \mc{A}_1(x_1,\ldots ,x_n;$
$\omega),\ldots ,\mc{A}_f^\omega =\mc{A}_f(x_1,\ldots
,x_n;\omega)$ abelian $*$-subalgebras of $\mc{M}$ $\exists
\mc{Z}^\omega =\mc{Z}(x_1,$ $\ldots ,x_n;\omega)\subset \mc{M}$
containing $p$ self-adjoint elements, such that
$$\dist_{||\cdot ||_2}\left(x_j,\sum_{\overset{1\leq j_1,\ldots ,j_{t+1}\leq f}
{1\leq t\leq
d}}\mc{A}_{j_1}^\omega\mc{Z}^\omega\mc{A}_{j_2}^\omega\mc{Z}^\omega\ldots
\mc{A}_{j_t}^\omega\mc{Z}^\omega\mc{A}_{j_{t+1}}^\omega\right)<\omega\,\,\forall
1\leq j\leq n.$$
\end{defn}
Proposition \ref{prop5} gives an estimate of the free entropy of a
(finite) system of generators of a $II_1$-factor $\mc{M}$ which
can be asymptotically decomposed as
$$\lim_{\omega\rightarrow 0}\,^{||\cdot ||_2}\sum_{\overset{1\leq j_1,\ldots , j_{t+1}\leq
f}{ 1\leq t\leq d}} \mc{A}_{j_1}^\omega \mc{Z}^\omega
\mc{A}_{j_2}^\omega \mc{Z}^\omega\ldots \mc{A}_{j_t}^\omega
\mc{Z}^\omega \mc{A}_{j_{t+1}}^\omega.$$ As in the statement of
Proposition \ref{prop3}, the approximations in the $||\cdot
||_2$-norm (\ref{ap}) hold for every $\omega >0$ if the
$II_1$-factor can be decomposed as above.
\begin{prop}\label{prop5}
Let $z_1,\ldots ,z_p$ be self-adjoint elements of a
$\mbox{\!I\!I}_1$-factor $\mc{M}$ and let $(\mc{A}_v)_{1\leq v\leq
f}$ be a family of abelian subalgebras of $\mc{M}$. Let
$x_1,\ldots ,x_n$ be self-adjoint generators of $\mc{M}$ and
assume that there exist projections
$p^{(v)}_1,\ldots,p^{(v)}_{r_v}\in \mc{A}_v$ and complex
noncommutative polynomials $(\phi_j)_{1\leq j\leq n}$ of degrees
$\leq d$ (where $d\geq 1$ is fixed) in the variables $(z_u)_{1\leq
u\leq p}$ such that
\begin{equation}\label{ap}
\left|\left|x_j-\phi_j\left((p^{(v)}_i)_{\overset{1\leq i\leq
r_v}{ 1\leq v\leq f}},(z_u)_{1\leq u\leq p}\right)
\right|\right|_2<\omega,\,\,j=1,\ldots,n,
\end{equation}
where $\omega\in (0,a]$ is a given positive number, and such that
in all monomials of every $\phi_j$ the projections $p^{(v)}_i$ and
$p^{(w)}_k$ are separated by some $z_u$ if $v\ne w$. Then
\begin{equation}
\chi(x_1,\ldots,x_n)\leq C(n,p,a,d,f)+(n-p-2f)\log\omega ,
\end{equation}
where $a=\max\lbrace ||x_j||_2+1|1\leq j\leq n\rbrace$ and
$C(n,p,a,d,f)$ is a constant that depends only on $n,p,a,d,f$.
\end{prop}
\begin{proof} As in the proof of Proposition \ref{prop3} we can assume that
$\phi_j=\phi_j^*\,\forall 1\leq j\leq n$ and fix $R>0$. Consider
an arbitrary element
$$\left((B_j)_{1\leq j\leq n},(P^{(v)}_i)_{\overset{1\leq i\leq r_v}{1\leq v\leq f}},
(Z_u)_{1\leq u\leq p}\right)$$ of
$$\Gamma_R\left((x_j)_{1\leq j\leq n},(p^{(v)}_i)_{\overset{1\leq i\leq r_v}
{1\leq v\leq f}},(z_u)_{1\leq u\leq p};m,k,\epsilon\right)$$ for
some large integers $m$, $k$ and small $\epsilon >0$. Eventually
after further restricting $m$, $\epsilon$, we can find mutually
orthogonal projections $Q^{(v)}_1,\ldots ,Q^{(v)}_{r_v}\in
\mc{M}_k^{sa}$ with $\mbox{rank} (Q^{(v)}_i)=[\tau (p^{(v)}_i)k]$
$\forall 1\leq i\leq r_v$, such that
$$\left|\left|B_j-\phi_j\left((Q^{(v)}_i)_{\overset{1\leq i\leq r_v}{ 1\leq v\leq f}},
(Z_u)_{1\leq u\leq p}\right) \right|\right|_2<\omega\,\,\forall
1\leq j\leq n\,\,.$$ If $S^{(v)}_1,\ldots ,S^{(v)}_r\in
\mc{M}_k^{sa}$ are fixed, mutually orthogonal projections with
$\mbox{rank}(S^{(v)}_i)=[\tau (p^{(v)}_i)k]$ for every $1\leq
i\leq r_v$, then there exists a unitary $U^{(v)}\in \mc{U}(k)$
such that $Q^{(v)}_i=U^{(v)*}S_iU^{(v)}$ $\forall 1\leq i\leq
r_v$. The previous inequality becomes
$$\left|\left|B_j-\phi_j\left((S^{(v)}_i)_{\overset{1\leq i\leq r_v}{ 1\leq v\leq f}},
(Z_u)_{1\leq u\leq
p},(\mbox{Re}(U^{(v)}),\mbox{Im}(U^{(v)}))_{1\leq v\leq f}\right)
\right|\right|_2<\omega ,$$ and all the components of $\phi_j$ are
polynomials of degrees $\leq 3d+2$ in the last $p+2f$ variables.
Reasoning as in the last part of the proof of Proposition
\ref{prop3} we can easily obtain now the estimate
$\chi(x_1,\ldots,x_n)\leq C(n,p,a,d,f)+(n-p-2f)\log\omega$.
\end{proof}
\subsection{Abelian dimension of free group factors}
\begin{thm}\label{prop6}
If $n\geq p+2f+1$, then the free group factor $\mc{L}(\mb{F}_n)$
does not admit an asymptotic decomposition of the form
$$\lim_{\omega\rightarrow 0}\,^{||\cdot ||_2} \sum_{\overset{1\leq
j_1,\ldots , j_{t+1}\leq f}{1\leq t\leq d}} \mc{A}_{j_1}^\omega
\mc{Z}^\omega \mc{A}_{j_2}^\omega \mc{Z}^\omega \ldots
\mc{A}_{j_t}^\omega \mc{Z}^\omega \mc{A}_{j_{t+1}}^\omega ,$$
where each subset $\mc{Z}^\omega$ contains $p$ self-adjoint
elements, $\mc{A}_1^\omega,\ldots ,\mc{A}_f^\omega\subset
\mc{L}(\mb{F}_n)$ are abelian $*$-subalgebras and $d\geq 1$ is an
integer.
\end{thm}
\begin{proof} Apply Proposition \ref{prop5} in the same manner Proposition \ref{prop3}
was used in the proof of Theorem \ref{prop4}.
\end{proof}
\begin{cor}\label{cor9}
If $\mc{P}\subset \mc{L}(\mb{F}_n)$ is a subfactor of finite index
and if $n\geq p+2f+2$, then $\mc{P}$ can not be asymptotically
decomposed as
$$\lim_{\omega\rightarrow 0}\,^{||\cdot ||_2}\sum_{\overset{1\leq
j_1,\ldots , j_{t+1}\leq f}{ 1\leq t\leq d}} \mc{A}_{j_1}^\omega
\mc{Z}^\omega \mc{A}_{j_2}^\omega \mc{Z}^\omega \ldots
\mc{A}_{j_t}^\omega \mc{Z}^\omega \mc{A}_{j_{t+1}}^\omega ,$$
where each subset $\mc{Z}^\omega$ contains $p$ self-adjoint
elements of $\mc{P}$, $\mc{A}_1^\omega,\ldots ,\mc{A}_f^\omega$
$\subset \mc{P}$ are abelian $*$-subalgebras, and  $d\geq 1$ is an
integer.
\end{cor}
\begin{proof} It is a direct consequence of Theorem \ref{prop6} and of
decomposition $\mc{L}(\mb{F}_n)=\mc{P}e_\mc{Q}\mc{P}$ (see the
proof of Corollary \ref{cor6}).
\end{proof}
\begin{cor}\label{can}
If $n\geq p+2f+1$, then the free group factor $\mc{L}(\mb{F}_n)$
can not be decomposed as
$$\overline{\s}^w\sum_{\overset{1\leq j_1,\ldots ,j_{t+1}\leq
f}{ 1\leq t\leq d}} \mc{A}_{j_1}\mc{Z}\mc{A}_{j_2}\mc{Z}\ldots
\mc{A}_{j_t}\mc{Z}\mc{A}_{j_{t+1}} ,$$  where $\mc{Z}\subset
\mc{L}(\mb{F}_n)$ contains $p$ self-adjoint elements,
$\mc{A}_1,\ldots ,\mc{A}_f$ are abelian $*$-subalgebras of
$\mc{L}(\mb{F}_n)$, and $d\geq 1$ is an integer. Moreover, if
$\mc{P}\subset \mc{L}(\mb{F}_n)$ is a subfactor of finite index
and if $n\geq p+2f+2$, then $\mc{P}$ also can not be decomposed as
$$\overline{\s}^w\sum_{\overset{1\leq j_1,\ldots ,j_{t+1}\leq
f}{1\leq t\leq d}} \mc{A}_{j_1} \mc{Z}\mc{A}_{j_2}\mc{Z} \ldots
\mc{A}_{j_t} \mc{Z} \mc{A}_{j_{t+1}} ,$$ for any subset $\mc{Z}$
containing $p$ self-adjoint elements of $\mc{P}$, any
$\mc{A}_1,\ldots ,$ $\mc{A}_f$ abelian $*$-subalgebras of
$\mc{P}$, and any integer $d\geq 1$.
\end{cor}
\begin{proof} Apply Theorem \ref{prop6} and Corollary \ref{cor9}, for
$\mc{Z}^\omega =\mc{Z}$, $\mc{A}_1^\omega =\mc{A}_1,\ldots ,$
$\mc{A}_f^\omega =\mc{A}_f$.
\end{proof}
Corollary \ref{cor4} settles the second part of the conjecture of
L. Ge and S. Popa (\cite{5}), in the case $n=\infty$. As a
reminder, $\ell_a(\mc{M})$ is defined as
$\min\{f\in\mb{N}\,\,|\,\,\exists \mc{A}_1,\ldots ,\mc{A}_f
\subset\mc{M}\,\,\mbox{abelian $*$-algebras s.t.}\,\,
\overline{\mbox{sp}}^w \mc{A}_1\mc{A}_2\ldots \mc{A}_f$
$=\mc{M}\}$ for every type $\mbox{\!I\!I}_1$-factor $\mc{M}$.
\begin{cor}\label{cor4}
$\ell_a(\mc{L}(\mb{F}_n))\geq [\frac{n-2}{2}]+1\,\,\forall\,4\leq
n\leq\infty$.
\end{cor}
\begin{proof} It follows from the first part of Corollary \ref{can}, for $\mc{Z}=\{1\}$.
\end{proof}
\begin{cor}\label{cor10}
If $\mc{P}\subset \mc{L}(\mb{F}_n)$ is a subfactor of finite index
and $5\leq n\leq\infty$, then $\ell_a(\mc{P})\geq
[\frac{n-3}{2}]+1$.
\end{cor}
\begin{proof} Apply the second part of Corollary \ref{can}.
\end{proof}
\begin{rem} One can combine both indecomposability properties
of $\mc{L}(\mb{F}_n)$ into a single statement: if $n\geq p+2f+1$,
then the free group factor $\mc{L}(\mb{F}_n)$ does not admit an
asymptotic decomposition of the form $$\lim_{\omega\rightarrow
0}\,^{||\cdot ||_2}\sum_{\overset{1\leq j_1,\ldots , j_{t+1}\leq
f}{1\leq t\leq d}} \mc{M}_{j_1}^\omega \mc{Z}^\omega
\mc{M}_{j_2}^\omega \mc{Z}^\omega \ldots \mc{M}_{j_t}^\omega
\mc{Z}^\omega \mc{M}_{j_{t+1}}^\omega ,$$ where each subset
$\mc{Z}^\omega$ contains $p$ self-adjoint elements, each
$\mc{M}_1^\omega,\ldots ,$ $\mc{M}_f^\omega\subset
\mc{L}(\mb{F}_n)$ is either a nonprime subfactor or an abelian
$*$-subalgebra and $d\geq 1$ is an integer.
\end{rem}
{\bf Acknowledgment.} The author would like to thank F. R\u
adulescu for suggestions and many helpful conversations.


\begin{thebibliography}{12}
\bibitem[Dy1]{1} Dykema, K.: {\em Interpolated free group factors}.
Pac. J. Math. 163 (1994), 123-135
\bibitem[Dy2]{2} Dykema, K.: {\em Two applications of free entropy}.
Math. Ann. 308 (1997), 547-558
\bibitem[Ge1]{18} Ge, L.: {\em Applications of free entropy to finite
von Neumann algebras}. Amer. J. Math. 119 (1997), 467-485
\bibitem[Ge2]{4} Ge, L.: {\em Applications of free entropy to finite
von Neumann algebras, $\mbox{\!I\!I}$}. Ann. of Math. (2) 147
(1998), 143-157
\bibitem[GePo]{5} Ge, L., Popa, S.: {\em On some decomposition properties for
factors of type $\mbox{\!I\!I}_1$}. Duke Math. J. 94 (1998),
79-101
\bibitem[Ha]{54} Haagerup, U.: {\em An Example of a Non Nuclear C$^*$-algebra
which has the Metric Approximation Property}. Invent. Math. 50
(1979), 279-293
\bibitem[Jo]{7} Jones, V. F. R.: {\em Index for Subfactors}. Invent. Math.
72 (1983), 1-25
\bibitem[JoSu]{8} Jones, V. F. R., Sunder, V. S.: {\em Introduction to
 subfactors}. New York, Cambridge University Press, 1997
\bibitem[Ka]{37} Kadison, R. V.: {\em Problems on von Neumann algebras}.
Baton Rouge Conference (1967), unpublished
\bibitem[KR]{50} Kadison, R. V., Ringrose, J.:
{\em Fundamentals of the Theory of Operator Algebras}, Vols. 1, 2. Academic
Press, Orlando, 1983, 1986.
\bibitem[MvN]{32} Murray, F., von Neumann, J.: {\em On rings of operators, $\mbox{\!I\!V}$}.
Ann. of Math. 44 (1943), 716-808
\bibitem[Po1]{53} Popa, S.: {\em Singular maximal abelian $*$-subalgebras in
continuous von Neumann algebras}. J. Funct. Analysis 50 (1983),
151-166
\bibitem[Po2]{55} Popa, S.: {\em Notes on Cartan subalgebras in type
$\mbox{\!I\!I}_1$ factors}. Math. Scand. 57 (1985), 171-188
\bibitem[Po3]{19} Popa, S.: {\em Free-independent sequences in type
$\mbox{\!I\!I}_1$ factors and related problems}. Ast\' erisque 232
(1995), 187-202
\bibitem[R\u a1]{36} R\u adulescu, F.: {\em The fundamental group of
$\mc{L}(\mb{F}_\infty)$ is $\mb{R}_+\setminus\{0\}$}. J. Am. Math.
Soc. 5 (1992), 517-532
\bibitem[R\u a2]{10} R\u adulescu, F.: {\em Random matrices, amalgamated free
products and subfactors of the von Neumann algebra of a free
group, of noninteger index}. Invent. Math. 115 (1994), 347-389
\bibitem[StZs]{40} Str\u atil\u a, \c S., Zsid\'{o}, L.: {\em Lectures on von
Neumann Algebras}. Editura Academiei, Bucure\c sti, Rom\^ ania, and
Abacus Press, Tunbridge Wells, Kent, England, 1979
\bibitem[Sz]{12} Szarek, S. J.: {\em Nets of Grassmann manifolds and orthogonal
group}. Proceedings of Research Workshop on Banach Space Theory (Bor-Luh-Lin,
ed.), The University of Iowa, June 29-31 (1981), 169-185
\bibitem[\c St]{11} \c Stefan, M. B.: {\em The primality of subfactors of
finite index in the interpolated free group factors}. Proc. of the
AMS 126 (1998), 2299-2307
\bibitem[Vo1]{13} Voiculescu, D.: {\em Circular and
semicircular systems and free  product factors}. Operator
Algebras, Unitary Representations, Enveloping  Algebras, and
Invariant Theory, Progress in Mathematics, Volume 92,
Birkh\"{a}user, Boston (1990), 45-60
\bibitem[Vo2]{14} Voiculescu, D.: {\em The analogues of entropy and of
Fisher's information measure in free probability theory,
$\mbox{\!I\!I}$}. Invent. Math. 118 (1994), 411-440
\bibitem[Vo3]{15} Voiculescu, D.: {\em The analogues of entropy and of
Fisher's information measure in free probability theory,
$\mbox{\!I\!I\!I}$: the absence of Cartan subalgebras}. G.A.F.A.
Vol. 6, No. 1 (1996), 172-199
\bibitem[Vo4]{21} Voiculescu, D.: {\em The analogues of entropy and of
Fisher's information measure in free probability theory,
$\mbox{\!I\!V}$: maximum entropy and freeness}. Free Probability
Theory (D. V. Voiculescu, ed.), Fields Institute Communications 12
(1997), 293-302
\bibitem[Vo5]{17} Voiculescu, D.: {\em A Strengthened Asymptotic Freeness Result
for Random Matrices with Applications to Free Entropy}. IMRN No. 1
(1998), 41-63
\bibitem[VDN]{16}Voiculescu, D. V., Dykema, K. J., Nica, A.: {\em Free Random
Variables}. CRM Monograph Series, AMS 1992
\bibitem[We]{22}Weyl, H.: {\em On the Volume of Tubes}. Amer. J. Math.
61 (1939), 461-472
\bibitem[vdW]{20} Waerden, B. L. van der: {\em Modern Algebra}, Vol.
\textbf{2}. New York, F. Ungar Pub. Co. 1949-1950
\end{thebibliography}
\end{document}